\documentclass[preprint,12pt]{article}

\usepackage{graphicx}

\usepackage{amssymb, amsmath}
\usepackage{amsthm}
\usepackage{bbm}
\usepackage{lmodern}

\newcommand{\R}{\mathbb{R}}
\newcommand{\N}{\mathbb{N}}

\newcommand{\PP}{\mathbf{P}}
\newcommand{\E}{\mathbf{E}}
\newcommand{\1}{\mathbbm{1}}
\newcommand{\F}{\mathfrak{F}}
\newcommand{\Fo}{\mathfrak{F}^{Y}}
\newcommand{\M}{\mathcal{M}}
\newcommand{\PR}{\mathcal{M}_{1}}
\newcommand{\demo}{\paragraph{Proof}}
\newcommand{\findemo}{\hfill$\Box$\\}

\newcommand{\te}{t^{*}}
\newcommand{\ts}{t^{*}}
\def\union{\mathop{\cup}}
\def\inter{\mathop{\cap}}

\newtheorem{theorem}{Theorem}[section]
\newtheorem{lemme}[theorem]{Lemma}
\newtheorem{proposition}[theorem]{Proposition}
\newtheorem{definition}[theorem]{Definition}
\newtheorem{hyp}[theorem]{Assumption}
\newtheorem{remarque}[theorem]{Remark}
\newtheorem{corollaire}[theorem]{Corollary}
\usepackage{vmargin}

\begin{document}




\title{Optimal stopping for partially observed piecewise-deterministic Markov processes\footnote{This work was supported by ARPEGE program of the French National Agency of Research (ANR), project FAUTOCOES, number ANR-09-SEGI-004.}}

\author{Adrien Brandejsky\and  Beno\^\i te de Saporta \and Fran\c cois Dufour}
\date{}
\maketitle

\begin{abstract}
This paper deals with the optimal stopping problem under partial observation for piecewise-deterministic Markov processes. We first obtain a recursive formulation of the optimal filter process and derive the dynamic programming equation of the partially observed optimal stopping problem. Then, we propose a numerical method, based on the quantization of the discrete-time filter process and the inter-jump times, to approximate the value function and to compute an actual $\epsilon$-optimal stopping time. We prove the convergence of the algorithms and bound the rates of convergence. 
\end{abstract}


60G40, 60J25, 93E20, 93E25, 93E10, 60K10




\section{Introduction}
The aim of this paper is to investigate an optimal stopping problem under partial observation for piecewise-deterministic Markov processes (PDMP) both from the theoretical and numerical points of view.
PDMP's have been introduced by Davis \cite{davis93} as a general class of stochastic models. They form a family of Markov processes involving deterministic motion punctuated by random jumps. The motion depends on three local characteristics, the flow $\Phi$,  the jump rate $\lambda$ and the transition measure $Q$, which selects the post-jump location. Starting from the point $x$, the motion of the process $(X_{t})_{t\geq 0}$ follows the flow $\Phi(x,t)$ until the first jump time $T_{1}$, which occurs either spontaneously in a Poisson-like fashion with rate $\lambda(\Phi(x,t))$ or when the flow hits the boundary of the state space. In either case, the location of the process at $T_{1}$ is selected by the transition measure $Q(\Phi(x,T_{1}),\cdot)$ and the motion restarts from $X_{T_{1}}$. We define similarly the time until the next jump and the next post-jump location and so on.
One important property of a PDMP, relevant for the approach developed in this paper, is that its distribution is completely characterized by the discrete time Markov chain
$(Z_{n},S_{n})_{n\in\N}$ where $Z_{n}$ is the $n$-th post-jump location and $S_{n}$ is the $n$-th inter-jump time.
A suitable choice of the state space and local characteristics provides stochastic models covering a large number of applications such as operations research \cite[section 33]{davis93}, reliability \cite{JRR12}, neurosciences \cite{pakdaman10}, internet traffic \cite{chafai10},
finance \cite{bauerle11}. This list of examples and references is of course not exhaustive.

In this paper, we consider an optimal stopping problem for a partially observed PDMP $(X_{t})_{t\geq 0}$.
Roughly speaking, the observation process $(Y_{t})_{t\geq 0}$ is a point process defined through the embedded discrete time Markov chain
$(Z_{n},S_{n})_{n\in\N}$. The inter-arrival times are given by $(S_{n})_{n\in\N}$ and the marks by  a noisy function of $(Z_{n})_{n\in\N}$.
For a given reward function $g$ and a computation horizon $N\in\N$, we study the following optimal stopping problem 
$$\sup_{\sigma\leq T_{N}}\E\left[g(X_{\sigma})\right],$$
where $T_N$ is the $N$-th jump time of the PDMP $(X_{t})_{t\geq 0}$, $\sigma$ is a stopping time with respect to the natural filtration $\Fo=(\Fo_{t})_{t\geq 0}$ generated by the observations $(Y_{t})_{t\geq 0}$.
In some applications, it may be more appropriate to consider a fixed optimization horizon $t_f$ rather than the random horizon $T_N$. This is a difficult problem with few references in the literature, see for instance \cite{edoli10} where the underlying process is not piecewise deterministic. Regarding PDMP's, this problem could be addressed using the same ideas as in \cite{brandejsky11}. It involves the time-augmented process $(X_t,t)$. Although this process is still a PDMP, its local characteristics may not have the same good properties as those of the original process leading to several new technical difficulties.

A general methodology to solve such a problem is to split it into two sub-problems.
The first one consists in deriving the filter process given by the conditional expectation of $X_{t}$ with respect to the observed information $\Fo_{t}$.
Its main objective is to transform the initial problem into a completely observed optimal stopping problem where the new state variable is the filter process.
The second step consists in solving this reformulated problem, the new difficulty being its infinite dimension.
Indeed, the filter process takes values in a set of probability measures.


Our work is inspired by~\cite{pham05} which deals with an optimal stopping problem under partial observation for a Markov chain with finite state space. The authors study the optimal filtering and 
convert their original problem into a standard optimal stopping problem for a continuous state space Markov chain.
Then they propose a discretization method based on a quantization technique to approximate the value function.
However, their method cannot be directly applied to our problem for the following main reasons related to
the specificities of PDMPs.

Firstly, PDMPs are continuous time processes. 
Although the dynamics can be described by the discrete-time Markov chain $(Z_{n},S_{n})_{n\in\N}$, this optimization problem remains intrinsically a \textit{continuous-time} optimization problem. Indeed, the performance criterion is maximized over the set of stopping times defined with respect to the \textit{continuous-time} filtration $(\Fo_{t})_{t\geq 0}$.
Consequently, our problem cannot be converted into a fully discrete time problem.

Secondly, the distribution of a PDMP combines both absolutely continuous and singular components. This is due to the existence of forced jumps when the process hits the boundary of the state space. As a consequence the derivation of the filter process is not straightforward. In particular, the absolute continuity hypothesis \textbf{(H)} of \cite{pham05} does not hold.


Thirdly, in our context the reformulated optimization problem is not standard, unlike in \cite{pham05}.
As already explained, this reformulated optimization problem
 combines \textit{continuous-time} and \textit{discrete-time} features.
Consequently, this problem does not correspond to the classical optimal stopping problem of a discrete-time Markov chain.
Moreover, it is different from the optimal stopping problem of a PDMP under complete observation mainly because the new state variables given by the Markov chain $(\Pi_{n},S_{n})_{n\geq 0}$ are not the underlying Markov chain of some PDMP. Therefore the results of the literature \cite{saporta10,gugerli86,pham05} cannot be used.

%

Finally, 
a natural way to proceed with the numerical approximation is then to follow the ideas developed in \cite{saporta10,pham05} namely to replace the filter $\Pi_{n}$ and the inter-jump time $S_{n}$ by some finite state space approximations in the dynamic programming equation. However, a noticeable difference from \cite{saporta10} lies in the fact that the dynamic programming operators therein were Lipschitz continuous whereas our new operators are only Lipschitz continuous between some points of discontinuity. We overcome this drawback by splitting the operators into their restrictions onto their continuity sets.
This way, we obtain not only an approximation of the value function of the optimal stopping problem but also an $\epsilon$-optimal stopping time with respect to
the filtration $(\Fo_{t})_{t\geq 0}$ that can be computed in practice.

Our approximation procedure for random variables is based on quantization. There exists an extensive literature on this method. The interested reader may for instance consult \cite{gray98,pages04} and the references within. The quantization of a random variable $X$ consists in finding a finite grid such that the projection $\widehat{X}$ of $X$ on this grid minimizes some $L^{p}$ norm of the difference $X-\widehat{X}$. Roughly speaking, such a grid will have more points in the areas of high density of $X$. As explained for instance in \cite[section 3]{pages04}, under some Lipschitz-continuity conditions, bounds for the rate of convergence of functionals of the quantized process towards the original process are available, which makes this technique especially appealing.  Quantization methods have been developed recently in numerical probability or optimal stochastic control with applications in finance, see e.g. \cite{pages04,bally03,bally05}. 

The paper is organized as follows. Section \ref{section-def} introduces the notation, recalls the definition of a PDMP, presents our assumptions and defines the optimal stopping problem we are interested in, especially the observation process. The recursive formulation of the filter process is derived in Section~\ref{section-filtre}. In Section \ref{section-dynamic}, we reduce our partially observed problem for the PDMP $(X_{t})_{t\geq 0}$ to a completely observed one involving the process $(\Pi_{n},S_{n})_{n\in\N}$ for which we provide the dynamic programming equation and construct a family of $\epsilon$-optimal stopping times. Then, our numerical methods to compute the value function and an $\epsilon$-optimal stopping time are presented in Section \ref{section-quantif} where we also prove the convergence of our algorithms after having recalled the main features of quantization. Finally, an academic example is discussed in Section \ref{section-example} while technical results are postponed to the Appendices.
%
\section{Definition and notation}\label{section-def}
In this first section, let us define a piecewise-deterministic Markov process (PDMP) and introduce some general assumptions.
For any metric space $E$, we denote $\mathcal{B}(E)$ its Borel $\sigma$-field, $B(E)$ the set of real-valued, bounded and measurable functions defined on~$E$ and $BL(E)$ the subset of functions of $B(E)$ that are Lipschitz continuous. For $a,b \in \R$, denote $a\wedge b=\min(a,b)$ and $a\vee b=\max(a,b)$.
\subsection{Definition of a Piecewise-Deterministic Markov Process}
Let $E$ be an open subset of $\R^{d}$. Let $\partial E$ be its boundary and $\overline{E}$ its closure and for any subset $A$ of $E$, $A^{c}$ denotes its complement. 
A PDMP is defined by its local characteristics $(\Phi,\lambda,Q)$.
\begin{itemize}
\item{The flow $\Phi : \mathbb{R}^{d}\times \mathbb{R}^{+}\rightarrow \mathbb{R}^{d}$ is continuous. For all $t\in \mathbb{R}^{+}$, $\Phi(\cdot,t)$ is an homeomorphism and $t\rightarrow \Phi(\cdot,t)$ is a semi-group:  for all $x \in \mathbb{R}^{d}$, $\Phi(x,t+s)=\Phi(\Phi(x,s),t)$.
    For all $x \in E$, define the deterministic exit time from $E$:
    $t^*(x)=
    \inf\{t>0 \text{ such that } \Phi(x,t)\in\partial E\}.$
    We use here and throughout the convention $\inf \emptyset = + \infty$.
    }
\item{The jump rate $\lambda : \overline{E}\rightarrow \mathbb{R}^+$ is measurable and satisfies:
$$\forall x \in E\text{,  }\exists \epsilon >0\text{ such that }\int_0^\epsilon \lambda(\Phi(x,t))dt< +\infty.$$
}
\item{Finally, $Q$ is a Markov kernel on $(\overline{E},\mathcal{B}(\overline{E}))$ which satisfies:
$$\forall x \in \overline{E}\text{,  }Q(x,E\backslash\{x\})=1.$$
}
\end{itemize}
From these characteristics, it can be shown \cite{davis93} that there exists a filtered probability space $(\Omega,\mathcal{F},(\mathcal{F}_t)_{t\in\R^{+}},(\mathbf{P}_x)_{x\in E})$ on which a process $(X_t)_{t\in \mathbb{R}^+}$ is defined. Its motion, starting from a point $x\in E$, may be constructed as follows.
Let $T_1$ be a nonnegative random variable with survival function:
$$\PP_x(T_1>t)=
\left\{\begin{array}{ll}
e^{-\Lambda(x,t)} & \text{if } 0\leq t<t^*(x), \\
0 & \text{if }t\geq t^*(x),
\end{array}\right.$$
where for $x\in E$ and $t\in [0,t^*(x)]$,
$\Lambda(x,t)=\int_0^t\lambda(\Phi(x,s))ds.$
One then chooses an $E$-valued random variable $Z_1$ with distribution $Q(\Phi(x,T_1),\cdot)$. The trajectory of $X_t$ for $t\leq T_1$ is:
$$X_t=\left\{
\begin{array}{ll}
\Phi(x,t)&\text{ if }t<T_1,\\
Z_1&\text{ if }t=T_1.
\end{array}\right.
$$
Starting from the point $X_{T_1}=Z_1$, one selects in a similar way $S_{2}=T_2-T_1$ the time between $T_{1}$ and the next jump time $T_2$, as well as $Z_2$ the next post-jump location and so on. Davis showed \cite{davis93} that the process so defined is a strong Markov process $(X_t)_{t\geq 0}$ with jump times $(T_n)_{n\in {\mathbb{N}}}$ ($T_0=0$). The process $(Z_n,S_n)_{n\in\N}$ where $Z_n=X_{T_n}$ is the $n$-th post-jump location and $S_n=T_n-T_{n-1}$ ($S_0=0$) is the $n$-th inter-jump time is clearly a discrete-time Markov chain.
\subsection{Notation and assumptions}
The following non explosion assumption about the jump-times is standard (see for example \cite[section 24]{davis93}).
\begin{hyp}\label{hyp-Tk_goes_to_infty}
For all $(x,t)\in E\times \mathbb{R}^+$, $\E_x\left[\sum_k \mathbbm{1}_{\{T_k<t\}}\right]<+\infty$.
\end{hyp}
It implies that $T_k\rightarrow +\infty$ a.s. when $k\rightarrow +\infty$.
Moreover, we make the following assumption about the transition kernel $Q$.
\begin{hyp}
We assume that there exists a finite set $E_{0}=\{x_{1},\ldots,x_{q}\}\subset E$ such that for all $x\in E$, one has $Q(x,E_{0})=1$.
\end{hyp}
In other words, for all $n\in\N$, $Z_{n}$ may only take its values in the finite set~$E_{0}$. This assumption ensures that the filter process, defined in the next section, has finite dimension. This is required to derive a tractable numerical method in Section~\ref{section-quantif}. When this assumption does not hold, one may consider a preliminary discretization of the transition kernel to introduce it.
\begin{hyp}\label{hyp-ts-bounded}
We assume that the function $\ts$ is bounded on $E_{0}$ i.e. for all $m\in \{1,\ldots,q\}$, we assume that $0<\ts(x_{m})<+\infty$.
\end{hyp}
\begin{definition} \label{def-ts-order}
For all $m\in \{1,\ldots,q\}$, denote $\ts_{m}=\ts(x_{m})$ and assume that $x_{1}$,\ldots, $x_{q}$ are numbered such that $\ts_{1} \leq \ts_{2}\leq$ \ldots $\leq \ts_{q}$. Moreover, let $\ts_{0}=0$.
\end{definition}
For any function $w$ in $B(E)$, introduce the following notation 
\begin{equation*}
Qw(x)=\int_E w(y)Q(x,dy)=\sum_{i=1}^{q}w(x_{i})Q(x,x_{i}),\qquad C_w=\sup_{x\in \overline{E}}|w(x)|.
\end{equation*}
For any Lipschitz continuous function $w$ in $BL(E)$, denote $[w]$ its Lipschitz constant 
$$[w]=\sup_{x\neq y \in E}\frac{|w(x)-w(y)|}{|x-y|}.$$
\begin{hyp}
The jump rate $\lambda$ is in $B(\overline{E})$ i.e. is bounded by $C_{\lambda}$.
\end{hyp}
Denote $\M(E_{0})$ the set of finite signed measures on $E_{0}$ and $\PR(E_{0})$ the subset of probability measures on $E_{0}$. We equip $\M(E_{0})$ with the norm $|\cdot|$ given by $|\pi|=\sum_{i=1}^{q}|\pi^{i}|$ where $\pi^{i}$ denotes $\pi(\{x_{i}\})$.
\subsection{Partially observed optimal stopping problem}
We consider from now on a PDMP $(X_{t})_{t\geq 0}$ of which the initial state $X_{0}=Z_{0}$ is a fixed point $x_{0}\in E_0$. We assume that this PDMP is observed through a noise and we now turn to the description of our observation procedure.
For all $n\in\N$, we assume that $S_{n}$ is perfectly observed but that $Z_{n}$ is not (except for the initial state $Z_{0}$). In some examples, it seems reasonable to consider that the jump times of the process are observed (for instance, if the jumps correspond to changes of environment) and that, when a jump occurs, the actual post-jump location is measured with a noise. The \textit{observation} process of $Z_{n}$, denoted by $Y_{n}$ is assumed to be of the following form: $Y_{0}=x_{0}$ (deterministic) and for $n\geq 1$,
\begin{equation}\label{def-Yn}
Y_{n}=\varphi(Z_{n})+W_{n},
\end{equation}
where $\varphi : E_{0} \rightarrow \R^{d}$ and where the \textit{noise} $(W_{n})_{n\geq 1}$ is a sequence of $\R^{d}$-valued, i.i.d. random variables with bounded density function $f_{W}$ that are also independent from $(Z_{n},S_{n})_{n\in\N}$. 
In order to define real-valued stopping times adapted to the observation process, we need to consider a continuous time version of the observation process. We therefore define the piecewise-constant process $(Y_t)_{t\geq 0}$ with a slight abuse of notation\footnote{The quantity $Y_{n}$ represents the value of the process $(Y_{t})_{t\geq 0}$ at time $t=T_{n}$ and must not be confused with the value of the process at time $t=n$.} as
$$Y_{t}=\sum_{j=0}^{+\infty}\1_{[T_{j},T_{j+1}[}(t)Y_{j}.$$ 
Let $\Fo=(\Fo_{t})_{t\geq 0}$ be the filtration generated by $(Y_{t})_{t\geq 0}$ (the \textit{observed} filtration) and $\F=(\F_{t})_{t\geq 0}$ be the filtration generated by $(X_{t},Y_{t})_{t\geq 0}$ (the \textit{total} filtration). Without changing the notation, we then complete these filtrations with all the $\PP$-null sets. This leads us to the following definition.
\begin{definition} Denote $\Sigma^{Y}$ the set of $(\Fo_{t})_{t\geq 0}$-stopping times that are a.s. finite and for $n\in\N$, define
$$\Sigma_{n}^{Y}=\left\{ \sigma\in \Sigma^{Y}  \text{ such that } \sigma\leq T_{n}  \text{ a.s.} \right\}.$$
\end{definition}
For all $n\in\N$, we define the filter $\Pi_{n}\in\PR(E_{0})$. The quantity $\Pi_{n}(\{x_{i}\})$, denoted by $\Pi_{n}^{i}$, represents the probability of the event $\{Z_{n}=x_{i}\}$ given the information available until time $T_{n}$ i.e.
\begin{equation}\label{def-filtre}
\forall i\in\{1,\ldots,q\},\qquad \Pi_{n}^{i}=\E[\1_{\{Z_{n}=x_{i}\}}\big| \Fo_{T_{n}}].
\end{equation}
Finally, let $N\in\N$ be the \textit{horizon} and $g\in B(\overline{E})$ the \textit{reward function}, we are interested in maximizing the following \textit{performance criterion}
$$\E\left[g(X_{\sigma})\big|\Pi_0=\pi\right]$$
with respect to the stopping times $\sigma\in \Sigma_{N}^{Y}$.
The \textit{value function} associated to this partially observed optimal stopping problem is given by
\begin{equation}\label{opt-stop-pb}
v(\pi)=\sup_{\sigma\in \Sigma_{N}^{Y}}\E\left[g(X_{\sigma})\big|\Pi_0=\pi\right],
\end{equation}
where $\pi$ is a probability measure in $\PR(E_{0})$. The solution of our problem is then obtained by setting $\pi=\delta_{x_{0}}$.
For some applications, it would be interesting to consider a more general form for the reward function such as an integral term also possibly depending on the observation process, see for instance \cite{ludkovski12}. However, this new setup would lead to several technical difficulties. In particular, the dynamic programming would be more complex. Thus the derivation of the error bounds for the numerical approximation would be possibly intractable.

We will also need the following assumption about the reward function $g$ associated with the optimal stopping problem.
\begin{hyp}\label{hyp-g-lip}
The function $g$ is in $B(\overline{E})$ i.e. bounded by $C_{g}$ and there exists $[g]_{2}\in \mathbb{R}^+$ such that for all $i\in \{1,\ldots,q\}$ and $t,u\in [0,\ts_{i}]$, one has:
$$|g(\Phi(x_{i},t))-g(\Phi(x_{i},u))|\leq [g]_{2}|t-u|.$$
\end{hyp}
Now, the aims of this paper are 
first {to explicit the filter process $(\Pi_{n})_{n\in\N}$ (Section~\ref{section-filtre});}
 second {to rewrite the partially observed optimal stopping problem \eqref{opt-stop-pb} as a totally observed one for a suitable Markov chain on $\PR(E_{0})\times \R^{+}$ (Section~\ref{sec-opt-stop-complete});}
 third {to derive a dynamic programming equation and construct a family of $\epsilon$-optimal stopping times (Section~\ref{section-dynamic-eq});}
 and finally {to propose a numerical method to compute an approximation of the value function and an $\epsilon$-optimal stopping time (Section~\ref{section-quantif}).}
As a starting point, we will derive, in the next section, a recursive construction of the optimal filter that is the key point of our approach.
%
\section{Optimal filtering}\label{section-filtre}
The goal of this section is to obtain a recursive formulation of the filter~$\Pi_{n}$. As far as we know, there is no result concerning the filter process for generic PDMPs. We may however refer to \cite{arjas92} for a recursive formulation of the filter for point processes, that can be seen as a sub-class of PDMP's.
For all $n\in\N$, we denote $\mathcal{G}_{n}=(Y_{0},S_{0},\ldots,Y_{n},S_{n})$. The continuous-time observation process $(Y_{t})_{t\geq 0}$ being a point process in the sense developed in~\cite{bremaud81}, one has $\Fo_{T_{n}}=\sigma(\mathcal{G}_{n})$ (see \cite[page 58, Theorem T2]{bremaud81}). Moreover, $\F_{T_{n}}=\sigma(Z_{0},\ldots,Z_{n})\vee \Fo_{T_{n}}$. Concerning the filter $\Pi_{n}$, first notice that, since it is an $\Fo_{T_{n}}$-measurable random variable, there exists for all $n\in\N$ a measurable function $\pi_{n}:(\R^{d}\times \R^{+})^{n+1}\rightarrow \PR(E_{0})$ such that $\Pi_{n}=\pi_{n}(\mathcal{G}_{n}).$
As in the case of the Kalman-Bucy filter, the iteration leading from $\Pi_{n-1}$ to $\Pi_{n}$ can be split into two steps : prediction and correction. For all $n\geq 1$, let $\mu_{n}^{-}$ be the conditional distribution of $(Z_{n},S_{n})$ given $\Fo_{T_{n-1}}$. Thus, $\mu_{n}^{-}$ is a transition kernel defined on $(\R^{d}\times \R^{+})^{n} \times \mathcal{B}(E_{0}\times \R^{+})$ for all $j\in \{1,\ldots,q\}$ and $\gamma_{n-1} \in (\R^{d}\times\R^{+})^{n}$ by
\begin{equation}\label{def-mu-}
\mu_{n}^{-}(\gamma_{n-1},\{x_{j}\},ds)=\PP(Z_{n}=x_{j},S_{n}\in ds|\mathcal{G}_{n-1}=\gamma_{n-1}).
\end{equation}
\begin{lemme}\label{lemme-filtre-zys} 
For all $\gamma_{n-1} \in (\R^{d}\times\R^{+})^{n}$, we have the following equality of probability measures on $E_{0}\times  \R^{d}\times \R^{+}$, for all $j\in \{1,\ldots,q\}$,
$$\PP(Z_{n}=x_{j},Y_{n}\in dy,S_{n}\in ds|\mathcal{G}_{n-1}=\gamma_{n-1})=\mu_{n}^{-}(\gamma_{n-1},\{x_{j}\},ds)f_{W}(y-\varphi(x_{j}))dy.$$
\end{lemme}
\demo Set $h$ in $B(E_{0}\times  \R^{d}\times \R^{+})$, using Eq.~\eqref{def-Yn} that defines $Y_{n}$, one has
\begin{eqnarray*}\lefteqn{\E\left[h(Z_{n},Y_{n},S_{n})\big|\mathcal{G}_{n-1}=\gamma_{n-1}\right]}\\
&=& 
\sum_{j=1}^{q}\int h(x_{j},\varphi(x_{j})+w,s)\PP(Z_{n}=x_{j},S_{n}\in ds,W_{n}\in dw|\mathcal{G}_{n-1}=\gamma_{n-1}).
\end{eqnarray*}
Moreover, $W_{n}$ is independent from $\sigma(Z_{n},S_{n})\vee\Fo_{T_{n-1}}=\sigma(Z_{n},S_{n},\mathcal{G}_{n-1})$ and admits the density function $f_{W}$.
Consequently, one easily obtains  the result
by using the change of variable $y=\varphi(x_{j})+w$.
\findemo\\

Integrating w.r.t. to the first variable in the previous lemma (i.e. summing w.r.t. $x_{j}$) yields the following result.
\begin{corollaire}\label{lemme-filtre-ys}
For all $\gamma_{n-1} \in (\R^{d}\times\R^{+})^{n}$, we have the following equality of probability measures on $ \R^{d}\times \R^{+}$,
$$\PP(Y_{n}\in dy,S_{n}\in ds|\mathcal{G}_{n-1}=\gamma_{n-1})=\left[\sum_{j=1}^{q}\mu_{n}^{-}(\gamma_{n-1},\{x_{j}\},ds)f_{W}(y-\varphi(x_{j}))\right]dy.$$
\end{corollaire}
\begin{lemme}\label{lemme-filtre-mu}
For all $n\geq 1$, $\gamma_{n-1} \in (\R^{d}\times\R^{+})^{n}$ and $j\in\{1,\ldots,q\}$, the distribution $\mu_{n}^{-}$, defined by Eq.~\eqref{def-mu-}, satisfies
\begin{eqnarray*}
\lefteqn{\mu_{n}^{-}(\gamma_{n-1},\{x_{j}\},ds)}\\
&=&\sum_{m=0}^{q-1}\1_{\{s\in]\ts_{m};\ts_{m+1}[\}}\left(\sum_{i=m+1}^{q}\pi_{n-1}^{i}(\gamma_{n-1})\lambda(\Phi(x_i,s))e^{-\Lambda(x_i,s)}Q(\Phi(x_i,s),x_j)\right)ds\\
&&+\sum_{m=1}^{q}\left(\pi_{n-1}^{m}(\gamma_{n-1})e^{-\Lambda(x_m,\ts_m)}Q(\Phi(x_m,\ts_m),x_j)\right)\delta_{\ts_m}(ds).
\end{eqnarray*}
\end{lemme}
\demo Let $h$ be a function of $B(E_{0}\times \R^{+})$. Since $\sigma(\mathcal{G}_{n-1})=\Fo_{T_{n-1}}\subset\F_{T_{n-1}}$, the law of iterated conditional expectations yields
$$\E\left[h(Z_{n},S_{n})\big|\mathcal{G}_{n-1}=\gamma_{n-1}\right]=\E\left[\E\left[h(Z_{n},S_{n})\big|\F_{T_{n-1}}\right]\big|\mathcal{G}_{n-1}=\gamma_{n-1}\right].$$
Besides, $\F_{T_{n-1}}=\sigma(Z_{0},S_{0},W_{0},\ldots,Z_{n-1},S_{n-1},W_{n-1})$ so that $$\E\left[h(Z_{n},S_{n})\big|\F_{T_{n-1}}\right]=\E\left[h(Z_{n},S_{n})\big|Z_{0},S_{0},\ldots,Z_{n-1},S_{n-1}\right],$$ by independence of the sequences $(W_{n})_{n\in\N}$ and $(Z_{n},S_{n})_{n\in\N}$. Now, we apply the Markov property of $(Z_{n},S_{n})_{n\in\N}$ 
and a well-known special feature of the transition kernel of the underlying Markov chain of a PDMP to obtain
$$\E\left[h(Z_{n},S_{n})\big|\F_{T_{n-1}}\right]=\E\left[h(Z_{n},S_{n})\big|Z_{n-1},S_{n-1}\right]=\E\left[h(Z_{n},S_{n})\big|Z_{n-1}\right].$$
Moreover, the transition kernel can be explicitly expressed in terms of the local characteristics of the PDMP, and this yields the next equations

\begin{eqnarray*}\lefteqn{\E[h(Z_{n},S_{n})|\mathcal{G}_{n-1}=\gamma_{n-1}]}\\
&=&\E\Big[\sum_{i=1}^{q}\1_{\{Z_{n-1}=x_{i}\}}\E[h(Z_{n},S_{n})|Z_{n-1}=x_{i}]\big|\mathcal{G}_{n-1}=\gamma_{n-1}\Big]\\
&=&\E\Big[\sum_{i=1}^{q}\1_{\{Z_{n-1}=x_{i}\}}\sum_{j=1}^{q}\Big[\int_{\R^{+}} h(x_{j},s)\lambda(\Phi(x_i,s))e^{-\Lambda(x_i,s)}\1_{\{s<\ts_i\}}Q(\Phi(x_i,s),x_j)ds\\
&&+h(x_{j},\ts_i)e^{-\Lambda(x_i,\ts_i)}Q(\Phi(x_i,\ts_i),x_j)\Big]\big|\mathcal{G}_{n-1}=\gamma_{n-1}\Big]\\
&=&\sum_{j=1}^{q}\Big(\int_{\R^{+}} h(x_{j},s)\sum_{i=1}^{q}\pi_{n-1}^{i}(\gamma_{n-1})\lambda(\Phi(x_i,s))e^{-\Lambda(x_i,s)}\1_{\{s<\ts_i\}}Q(\Phi(x_i,s),x_j)ds\\
&&+\sum_{i=1}^{q}h(x_{j},\ts_i)\pi_{n-1}^{i}(\gamma_{n-1})e^{-\Lambda(x_i,\ts_i)}Q(\Phi(x_i,\ts_i),x_j)\Big).
\end{eqnarray*}
This can be written equivalently as
\begin{eqnarray*}\lefteqn{\E\left[h(Z_{n},S_{n})\big|\mathcal{G}_{n-1}=\gamma_{n-1}\right]}\\
&=&\sum_{j=1}^{q}\bigg(\sum_{m=0}^{q-1}\bigg(\int_{\ts_{m}}^{\ts_{m+1}} h(x_{j},s)\sum_{i=m+1}^{q}\pi_{n-1}^{i}(\gamma_{n-1})\lambda(\Phi(x_i,s))e^{-\Lambda(x_i,s)}Q(\Phi(x_i,s),x_j)\bigg)ds\\
&&+\sum_{i=1}^{q}h(x_{j},\ts_i)\pi_{n-1}^{i}(\gamma_{n-1})e^{-\Lambda(x_i,\ts_i)}Q(\Phi(x_i,\ts_i),x_j)\bigg).
\end{eqnarray*}
Hence the result.
\findemo\\

We now state the main result of this section, namely the recursive formulation of the filter sequence $(\Pi_{n})_{n\in\N}$.
\begin{proposition}\label{prop-rec-filtre}
Let $\Psi=(\Psi^1,\ldots,\Psi^q) : \PR(E_{0})\times  \R^{d} \times \R^{+} \rightarrow \PR(E_{0})$ be defined as follows: for all $j\in\{1,\ldots,q\}$,
$$\Psi^{j}(\pi,y,s)=\sum_{m=0}^{q-1}\1_{\{s\in]\ts_{m};\ts_{m+1}[\}}\frac{\Psi^{j}_{m}(\pi,y,s)}{\overline{\Psi}_{m}(\pi,y,s)}+\sum_{m=1}^{q}\1_{\{s=\ts_{m}\}}\frac{\Psi^{*j}_{m}(y)}{\overline{\Psi}^{*}_{m}(y)},$$
where
\begin{eqnarray*}
\Psi^{j}_{m}(\pi,y,s)&=&\sum_{i=m+1}^{q}\pi^{i}\lambda(\Phi(x_i,s))e^{-\Lambda(x_i,s)}Q(\Phi(x_i,s),x_j)f_{W}(y-\varphi(x_{j})),\\
\overline{\Psi}_{m}(\pi,y,s)&=&\sum_{k=1}^{q}\Psi^{k}_{m}(\pi,y,s),\\
\Psi^{*j}_{m}(y)&=&Q(\Phi(x_m,\ts_m),x_j)f_{W}(y-\varphi(x_{j})),\\
\overline{\Psi}^{*}_{m}(y)&=&\sum_{k=1}^{q}\Psi^{*k}_{m}(y).
\end{eqnarray*}
Then, the filter, defined in Eq.~\eqref{def-filtre}, satisfies $\Pi_0^j=\PP(Z_0=x_j)$ and the following recursion: for all $n\geq 1$,
$$ \PP\text{-a.s.,}Pi_{n}=\Psi(\Pi_{n-1},Y_{n},S\qquad\_{n}).$$
\end{proposition}
\demo Fix $\gamma_{n-1}$ in $(\R^{d}\times\R^{+})^{n}$. Bayes formula yields for all $j\in\{1,\ldots,q\}$,
\begin{eqnarray*}
\lefteqn{\PP(Z_{n}=x_{j},Y_{n}\in dy,S_{n}\in ds\big|\mathcal{G}_{n-1}=\gamma_{n-1})=}\\
&&  \PP\big(Z_{n}=x_{j}\big|\mathcal{G}_{n}=(\gamma_{n-1},y,s)\big) \times  \PP(Y_{n}\in dy,S_{n}\in ds\big|\mathcal{G}_{n-1}=\gamma_{n-1}).
\end{eqnarray*}
Lemma \ref{lemme-filtre-zys} and Corollary \ref{lemme-filtre-ys} yield
\begin{multline*}\mu_{n}^{-}(\gamma_{n-1},\{x_{j}\},ds)f_{W}(y-\varphi(x_{j}))dy\\=\PP\left(Z_{n}=x_{j}\big|\mathcal{G}_{n}=(\gamma_{n-1},y,s)\right) \left[\sum_{k=1}^{q}\mu_{n}^{-}(\gamma_{n-1},\{x_{k}\},ds)f_{W}(y-\varphi(x_{k}))\right]dy.
\end{multline*}
With respect to $y$, one recognizes the equality of two absolutely continuous measures which implies the equality a.e. of the density functions. Thus, one has for almost all $y\in \R^{d}$ w.r.t. the Lebesgue measure,
\begin{eqnarray}\label{eq-egal-mesures}
\lefteqn{\mu_{n}^{-}(\gamma_{n-1},\{x_{j}\},ds)f_{W}(y-\varphi(x_{j}))}\\
&=&\PP\left(Z_{n}=x_{j}\big|\mathcal{G}_{n}=(\gamma_{n-1},y,s)\right) \left[\sum_{k=1}^{q}\mu_{n}^{-}(\gamma_{n-1},\{x_{k}\},ds)f_{W}(y-\varphi(x_{k}))\right].\nonumber
\end{eqnarray}
Eq.~\eqref{eq-egal-mesures} states the equality of two measures of the variable $s\in\R^{+}$ that contain both an absolutely continuous part and some weighted Dirac measures.
Denote $g_{1}(y,s)\nu_{1}(ds)$ (respectively $g_{2}(y,s)\nu_{2}(ds)$) the left-hand (resp. right-hand) side term of the previous equality. Eq.~\eqref{eq-egal-mesures} means that for all function $F\in B(\R^{+})$ and for almost all $y\in \R^{d}$ w.r.t. the Lebesgue measure, one has
\begin{equation}\label{integral-fonction-test}
\int F(s) g_{1}(y,s)\nu_{1}(ds)=\int F(s) g_{2}(y,s)\nu_{2}(ds),
\end{equation}
Recall that, from Lemma \ref{lemme-filtre-mu}, the distribution $\mu_{n}^{-}(\gamma_{n-1},\{x_{j}\},ds)$ has a density on the interval $]\ts_{m};\ts_{m+1}[$ denoted by $f_{m}(\gamma_{n-1},x_{j},s)$ and given by
$$f_{m}(\gamma_{n-1},x_{j},s)=\sum_{i=m+1}^{q}\pi_{n-1}^{i}(\gamma_{n-1})\lambda(\Phi(x_i,s))e^{-\Lambda(x_i,s)}Q(\Phi(x_i,s),x_j).$$
First, take $F(s)=H(s)\1_{\{s\in ]\ts_{m};\ts_{m+1}[\}}$ in equation~\eqref{integral-fonction-test} with $H\in B(\R^{+})$.
One has
from equation (\ref{eq-egal-mesures})
\begin{multline*}
\int_{\ts_{m}}^{\ts_{m+1}}H(s)f_{m}(\gamma_{n-1},x_{j},s)f_{W}(y-\varphi(x_{j}))ds\\
=\int_{\ts_{m}}^{\ts_{m+1}}H(s)\PP\left(Z_{n}=x_{j}\big|\mathcal{G}_{n}=(\gamma_{n-1},y,s)\right)\sum_{k=1}^{q}f_{m}(\gamma_{n-1},x_{k},s)f_{W}(y-\varphi(x_{k}))ds,
\end{multline*}
and thus on $]\ts_{m};\ts_{m+1}[$, almost surely w.r.t. the Lebesgue measure, one has
$$\PP\left(Z_{n}=x_{j}\big|\mathcal{G}_{n}=(\gamma_{n-1},y,s)\right)=\frac{f_{m}(\gamma_{n-1},x_{j},s)f_{W}(y-\varphi(x_{j}))}{\sum_{k=1}^{q}f_{m}(\gamma_{n-1},x_{k},s)f_{W}(y-\varphi(x_{k}))}.$$
Finally, for $m\in\{1,\ldots,q\}$, choosing $F(s)=\1_{\{s=\ts_{m}\}}$ in Eq.~\eqref{integral-fonction-test} yields the equality of the weights at the point $\ts_{m}$ thus, using Lemma \ref{lemme-filtre-mu},
\begin{eqnarray*}
\lefteqn{\PP\left(Z_{n}=x_{j}\big|\mathcal{G}_{n}=(\gamma_{n-1},y,\ts_{m})\right)}\\
&=&\frac{\pi_{n-1}^{m}(\gamma_{n-1})e^{-\Lambda(x_m,\ts_m)}Q(\Phi(x_m,\ts_m),x_j)f_{W}(y-\varphi(x_{j}))}{\sum_{k=1}^{q}\pi_{n-1}^{m}(\gamma_{n-1})e^{-\Lambda(x_m,\ts_m)}Q(\Phi(x_m,\ts_m),x_k)f_{W}(y-\varphi(x_{k}))}\\
&=&\frac{Q(\Phi(x_m,\ts_m),x_j)f_{W}(y-\varphi(x_{j}))}{\sum_{k=1}^{q}Q(\Phi(x_m,\ts_m),x_k)f_{W}(y-\varphi(x_{k}))}.\\
\end{eqnarray*}
Thus there exists two measurable sets $N_{y}\subset \R^{d}$ and $N_{s}\subset \R^{+}\backslash\{\ts_{1},\ldots,\ts_{q}\}$, negligible w.r.t. the Lebesgue measures on $\R^{d}$ and $\R$ respectively, such that for all $\gamma_{n-1}\in (\R^{d}\times\R^{+})^{n}$, $y\in \R^{d}\backslash N_{y}$, $s\in\R^{+}\backslash N_{s}$, one has
\begin{equation}\label{eq-filtre-esp-cond}
\pi_{n}(\gamma_{n-1},y,s)=\Psi(\pi_{n-1}(\gamma_{n-1}),y,s).
\end{equation}
On the one hand, one has $\PP(Y_{n}\in N_{y})\leq \sum_{j=1}^{q}\PP(\varphi(x_{j})+W_{n}\in N_{y})=0$ by absolute continuity of the distribution of $W_{n}$. On the other hand, $\PP(S_{n}\in N_{s})=0$ because the distribution of $S_{n}$ is absolutely continuous on $\R^{+}\backslash\{\ts_{1},\ldots,\ts_{q}\}$ and one has $N_{s}\inter \{\ts_{1},\ldots,\ts_{q}\}=\emptyset$. We therefore conclude from Eq.~\eqref{eq-filtre-esp-cond} that $\PP$-a.s., one has
$\pi_{n}(\mathcal{G}_{n-1},Y_{n},S_{n})=\Psi(\pi_{n-1}(\mathcal{G}_{n-1}),Y_{n},S_{n}).$
 The result follows since $\PP$-a.s., one has $\pi_{n}(\mathcal{G}_{n-1},Y_{n},S_{n})=\Pi_{n}$ and $\pi_{n-1}(\mathcal{G}_{n-1})=\Pi_{n-1}$.
\findemo\\

This proposition will play a crucial part in the sequel. On the one hand, this result will enable us to prove the Markov property of the sequence $(\Pi_{n},S_{n})_{n\geq 0}$ w.r.t. the observed filtration. On the other hand, the recursive formulation allows for simulation of the process $(\Pi_{n})_{n\geq 0}$ which is crucial to obtain numerical approximations. 
Finally, notice that the specific structure of the PDMP appears in the recursive formulation of the filter which contains both an absolutely continuous part and some weighted points.
%
\section{Dynamic programming}\label{section-dynamic}
The main objective of this section is to derive the dynamic programming equation for the value function of the partially observed optimal stopping problem \eqref{opt-stop-pb}. The proof of this result can be roughly speaking decomposed into two steps. The first point consists in converting the partially observed optimal stopping problem into an optimal stopping problem under complete observation where the state variables are described by the \textit{discrete-time} Markov chain
$(\Pi_{n},S_{n})_{n\geq 0}$ (see Section \ref{sec-opt-stop-complete}). It is important to remark that under this new formulation, the optimization problem remains intrinsically a \textit{continuous-time} optimization problem because the performance criterion is maximized over the set of stopping times with respect to the \textit{continuous-time} filtration $(\Fo_{t})_{t\geq 0}$.
We show in the second step (see Section \ref{section-dynamic-eq}) that the value function 
associated to the optimal stopping problem \eqref{opt-stop-pb} can be calculated by iterating a functional operator, labelled $L$ (see Definition \ref{def-G-H-I-J-L}).
As a by-product, we also provide a family of $\epsilon$-optimal stopping times.

We would like to emphasize that the results obtained in this section are not straightforward to obtain due to the specific structure of this optimization problem. Indeed, as already explained, it combines \textit{continuous-time} and \textit{discrete-time} features.
Consequently, this problem does not correspond to the classical optimal stopping problem of a discrete-time Markov chain.
Moreover, it is different from the optimal stopping problem of a PDMP under complete observation mainly because the new state variables given by the Markov chain $(\Pi_{n},S_{n})_{n\geq 0}$ are not the underlying Markov chain of some PDMP. Therefore the results of the literature \cite{saporta10,gugerli86} cannot be used.

These derivations require some technical results about the structure of the stopping times in $\Sigma^{Y}_{N}$.
For the sake of clarity in exposition, they are presented in the Appendix \ref{STOP-time}.
%
We start with a technical preliminary result required in the sequel, investigating the Markov property of the filter process.
\begin{proposition}\label{prop-markov}
The sequences $(\Pi_{n},Y_{n},S_{n})_{n\in \N}$, $(\Pi_{n},S_{n})_{n\in \N}$ and $(\Pi_{n})_{n\in \N}$ are $(\Fo_{T_{n}})_{n\in\N}$-Markov chains.
\end{proposition}
\demo
Let $h\in B(\PR(E_{0})\times \R^{d} \times \R^{+})$. The law of iterated conditional expectations yields
\begin{eqnarray*}
\E[h(\Pi_{n},Y_{n},S_{n})|\Fo_{T_{n-1}}]
=\E\big[\E[h(\Pi_{n},Y_{n},S_{n})| \F_{T_{n-1}}]\big|\Fo_{T_{n-1}}\big].
\end{eqnarray*}
From Proposition \ref{prop-rec-filtre} and Eq.~\eqref{def-Yn} which defines $Y_{n}$ one obtains
\begin{eqnarray*}
\lefteqn{\E[h(\Pi_{n},Y_{n},S_{n})|\F_{T_{n-1}}]}\\
&=&\E\big[h\big(\Psi(\Pi_{n-1},\varphi(Z_{n})+W_{n},S_{n}),\varphi(Z_{n})+W_{n},S_{n}\big)\big| \F_{T_{n-1}}\big]\\
&=&\sum_{j=1}^{q}\int h\big(\Psi(\Pi_{n-1},\varphi(x_{j})+w,s),\varphi(x_{j})+w,s\big)\\
&&\times \PP(Z_{n}=x_{j},W_{n}\in dw,S_{n}\in ds|\F_{T_{n-1}}).
\end{eqnarray*}
Yet, $W_{n}$ is independent from $\sigma(Z_{n},S_{n})\vee\F_{T_{n-1}}$ and admits the density function $f_{W}$. As in the proof of Lemma \ref{lemme-filtre-zys} one thus obtains
\begin{eqnarray*}
\lefteqn{\E[h(\Pi_{n},Y_{n},S_{n})|\F_{T_{n-1}}]}\\
&=&\sum_{j=1}^{q}\int h\big(\Psi(\Pi_{n-1},y,s),y,s\big)\PP(Z_{n}=x_{j},S_{n}\in ds|\F_{T_{n-1}})f_{W}(y-\varphi(x_{j}))dy.
\end{eqnarray*}
Besides, we have $\PP(Z_{n}=x_{j},S_{n}\in ds|\F_{T_{n-1}})=\PP(Z_{n}=x_{j},S_{n}\in ds|Z_{n-1})$ as in the proof of Lemma \ref{lemme-filtre-mu}, so that one has
\begin{eqnarray*}
\lefteqn{\E[h(\Pi_{n},Y_{n},S_{n})|\F_{T_{n-1}}]}\\
&=&\sum_{i=1}^{q}\1_{\{Z_{n-1}=x_{i}\}}\sum_{j=1}^{q}\int \Big(\int_{0}^{\ts_{i}} h\big(\Psi(\Pi_{n-1},y,s),y,s\big)\lambda\big(\Phi(x_i,s)\big)e^{-\Lambda(x_i,s)}Q\big(\Phi(x_{i},s),x_{j}\big)ds\\
&&+ h\big(\Psi(\Pi_{n-1},y,\ts_{i}),y,\ts_{i}\big)e^{-\Lambda(x_i,\ts_{i})}Q\big(\Phi(x_{i},\ts_{i}),x_{j}\big)\Big)f_{W}(y-\varphi(x_{j}))dy.
\end{eqnarray*}
Take now the conditional expectation w.r.t. $\Fo_{T_{n-1}}$, to obtain
\begin{eqnarray*}
\lefteqn{\E[h(\Pi_{n},Y_{n},S_{n})|\F_{T_{n-1}}]}\\
&=&\sum_{i=1}^{q}\Pi_{n-1}^{i}\sum_{j=1}^{q}\int \Big(\int_{0}^{\ts_{i}} h\big(\Psi(\Pi_{n-1},y,s),y,s\big)\lambda\big(\Phi(x_i,s)\big)e^{-\Lambda(x_i,s)}Q\big(\Phi(x_{i},s),x_{j}\Big)ds\\
&&+ h\big(\Psi(\Pi_{n-1},y,\ts_{i}),y,\ts_{i}\big)e^{-\Lambda(x_i,\ts_{i})}Q\big(\Phi(x_{i},\ts_{i}),x_{j}\big)\Big)f_{W}(y-\varphi(x_{j}))dy.
\end{eqnarray*}
Hence $\E[h(\Pi_{n},Y_{n},S_{n})|\Fo_{T_{n-1}}]$ is merely a function of $\Pi_{n-1}$ yielding the result for the three processes.
\findemo
%
\subsection{Optimal stopping problem under complete observation}\label{sec-opt-stop-complete}
In this section, we show how our optimal stopping problem under partial observation for the process $(X_{t})_{t\geq 0}$ can be converted into an optimal stopping problem under complete observation involving the Markov chain $(\Pi_{n},S_{n})_{0\leq n\leq N}$.
More precisely, for a fixed stopping time $\sigma\in \Sigma^{Y}_{N}$, we show in Proposition \ref{lemme value fonction} that the performance criterion $\E[g(X_{\sigma})|\Pi_0=\pi]$ can be expressed in terms of the discrete-time Markov chain $(\Pi_{n},S_{n})_{0\leq n\leq N}$.
We would like to emphasize the following important fact. Although the performance criterion can be written in terms of \textit{discrete-time} process, the optimization problem remains intrinsically a \textit{continuous-time} optimization problem. Indeed, the performance criterion is maximized over the set of stopping times with respect to the \textit{continuous-time} filtration
$(\Fo_{t})_{t\geq 0}$.

\begin{proposition}\label{lemme value fonction}
Let $\sigma\in \Sigma^{Y}$ and $n\geq 1$. For all $\pi \in \PR(E_{0})$ one has
\begin{eqnarray*} 
\lefteqn{\E[g(X_{\sigma\wedge T_{n}})|\Pi_0=\pi]}\\
&=&
\sum_{k=0}^{n-1}\sum_{i=1}^{q}\E[\1_{\{T_{k}\leq \sigma\}}\1_{\{R_{k}<\te_{i}\}}g\circ\Phi(x_{i},R_{k})e^{-\Lambda(x_{i},R_{k})}\Pi_{k}^{i}|\Pi_0=\pi]\\
&&+\sum_{i=1}^{q}\E[\1_{\{T_{n}\leq \sigma\}}g(x_{i})\Pi_{n}^{i}|\Pi_0=\pi],
\end{eqnarray*}
where $(R_{k})_{k\in \N}$ is the sequence of non negative random variables associated to $\sigma$ as introduced in Theorem \ref{theo-bremaud-adapted}.
\end{proposition}
\demo
We split $\E[g(X_{\sigma\wedge T_{n}})|\Pi_0=\pi]$ into several terms depending on the position of $\sigma$ w.r.t. the jump times $T_k$
\begin{eqnarray*}
{\E[g(X_{\sigma\wedge T_{n}})\big|\Pi_0=\pi]}
&=&\sum_{k=0}^{n-1}\sum_{i=1}^{q}\E[\1_{\{T_{k}\leq \sigma<T_{k+1}\}}\1_{\{Z_{k}=x_{i}\}}g\circ\Phi(x_{i},R_{k})|\Pi_0=\pi]\\
&&+\sum_{i=1}^{q}\E[\1_{\{T_{n}\leq \sigma\}}\1_{\{Z_{n}=x_{i}\}}g(x_{i})|\Pi_0=\pi].
\end{eqnarray*}
For notational convenience, consider
$$\left\{\begin{array}{cl}
A_{k,i}&=\1_{\{T_{k}\leq \sigma<T_{k+1}\}}\1_{\{Z_{k}=x_{i}\}}g\circ\Phi(x_{i},R_{k}),\\
B_{i}&=\1_{\{T_{n}\leq \sigma\}}\1_{\{Z_{n}=x_{i}\}}g(x_{i}).
\end{array}\right.$$
On the one hand, one has $\E[B_{i}|\Fo_{T_{n}}]=g(x_{i})\1_{\{T_{n}\leq \sigma\}}\Pi_{n}^{i}$ since $\{T_{n}\leq \sigma\}\in\Fo_{T_{n}}$ (see for instance \cite[p. 298, Theorem T7]{bremaud81}).
On the other hand, to compute $\E[A_{k,i}|\Fo_{T_{k}}]$, we use Lemma \ref{lemme-tech-event} to obtain
\begin{align*}
\E[A_{k,i}|\Fo_{T_{k}}]
&=\1_{\{T_{k}\leq \sigma\}}g\circ\Phi(x_{i},R_{k})\E[\1_{\{S_{k+1}>R_{k}\}}\1_{\{Z_{k}=x_{i}\}}|\Fo_{T_{k}}]\\
&=\1_{\{T_{k}\leq \sigma\}}g\circ\Phi(x_{i},R_{k})\E\big[\1_{\{Z_{k}=x_{i}\}}\E[\1_{\{S_{k+1}>R_{k}\}}|\F_{T_{k}}]\big|\Fo_{T_{k}}\big]\\
&=\1_{\{T_{k}\leq \sigma\}}g\circ\Phi(x_{i},R_{k})\E[\1_{\{Z_{k}=x_{i}\}}\1_{\{R_{k}<\te(Z_{k})\}}e^{-\Lambda(Z_{k},R_{k})}|\Fo_{T_{k}}]\\
&=\1_{\{T_{k}\leq \sigma\}}g\circ\Phi(x_{i},R_{k})\1_{\{R_{k}<\te_{i}\}}e^{-\Lambda(x_{i},R_{k})}\Pi_{k}^{i}.
\end{align*}
Details to obtain the third line in the above computations are provided by Lemma~\ref{lemme-tech-esp-cond}.
The result follows.
\findemo
\subsection{Dynamic programming equation}\label{section-dynamic-eq}
Based on the new formulation, the main objective of this section is to derive the backward dynamic programming equation.
It involves some operators introduced in Definition \ref{def-G-H-I-J-L}.
By iterating the operator, labelled $L$, we define a sequence of real valued functions $(v_{n})_{0\leq n\leq N}$ in Definition \ref{def-vn}.
Theorem \ref{value-fonction} establishes that $v_{n}$ is the value function of our partially observed optimal stopping problem with horizon $T_{N-n}$ and in particular that $v_{0}$ is the value function of problem defined in equation~\eqref{opt-stop-pb}.

Another important result of this section is given by Theorem \ref{theo-Sn-epsilon} which constructs a sequence of $\epsilon$-optimal stopping times.
%
\begin{definition}
\label{def-G-H-I-J-L}
The operators $G:B(\PR(E_{0})) \rightarrow B(\PR(E_{0}) \times \R^{+})$,
$H:B(E) \rightarrow B(\PR(E_{0})\times  \R^{+})$,
$J:B(\PR(E_{0}))\times B(E) \rightarrow B(\PR(E_{0})\times \R^{+})$, and $L:B(\PR(E_{0}))\times B(E) \rightarrow B(\PR(E_{0}))$
are defined for all $(v,h)\in B(\PR(E_{0}))\times B(E)$ and $(\pi,u)\in \PR(E_{0})\times \R^{+}$ by
\begin{eqnarray*}
Gv(\pi,u)&=&\E[v(\Pi_{1})\1_{\{S_{1}\leq u\}}| \Pi_0=\pi],\\
Hh(\pi,u)&=&\E\big[\sum_{i=1}^qh\circ\Phi(x_i,u)\Pi_0^i\1_{\{u<t^*_i\}}\1_{\{S_{1} > u\}}| \Pi_0=\pi\big],\\
J(v,h)(\pi,u)&=&Hh(\pi,u)+Gv(\pi,u),\\
L(v,h)(\pi)&=&\sup_{u\geq 0} J(v,h)(\pi,u).
\end{eqnarray*}
\end{definition}
\begin{definition}\label{def-vn} 
The sequence $(v_{n})_{0\leq n\leq N}$ of real-valued functions is defined on $\PR(E_{0})$ by
$$\left\{\begin{array}{rl}
v_{N}(\pi)&=\sum_{i=1}^{q}g(x_{i})\pi^{i},\\
v_{n-1}(\pi)&=L(v_{n},g)(\pi),\quad 1\leq n\leq N.
\end{array}\right.$$
\end{definition}
%

The following Theorem is the main result of this section showing that the operator $L$ is the dynamic programming operator associated to the initial optimization problem.
\begin{theorem}\label{value-fonction} 
For all $1\leq n\leq N$ and $\pi\in \PR(E_{0})$, one has 
$$\sup_{\sigma\in \Sigma_{n}^{Y}}\E[g(X_{\sigma})|\Pi_0=\pi] = v_{N-n}(\pi).$$
\end{theorem}
\demo
The proof of this result is based on Proposition \ref{theo-value-fonction} and Theorem \ref{theo-Sn-epsilon}.
Proposition \ref{theo-value-fonction} proves that $v_{N-n}$ is an upper bound for the value function of the problem with horizon $T_{n}$.
The reverse inequality is derived in Theorem \ref{theo-Sn-epsilon} by constructing a sequence of $\epsilon$-optimal stopping times.
\findemo\\

\begin{proposition}\label{theo-value-fonction} 
For all $1\leq n\leq N$ and $\pi\in \PR(E_{0})$, one has 
$$\sup_{\sigma\in \Sigma_{n}^{Y}}\E[g(X_{\sigma})|\Pi_0=\pi]\leq v_{N-n}(\pi).$$
\end{proposition}
\demo
Let $\sigma\in\Sigma^{Y}$. Consider $(R_{k})_{k\in \N}$ the sequence associated to $\sigma$ as introduced in Theorem \ref{theo-bremaud-adapted}.
We prove the theorem by induction on $n$. For $n=1$, Proposition~\ref{lemme value fonction} yields
\begin{eqnarray}
\E[g(X_{\sigma\wedge T_{1}})|\Pi_0=\pi]
&=&{\sum_{i=1}^{q}\E[\1_{\{R_{0}<\te_{i}\}}g\circ\Phi(x_{i},R_{0})e^{-\Lambda(x_{i},R_{0})}\Pi_{0}^{i}|\Pi_0=\pi]} \nonumber \\
&&+{\sum_{i=1}^{q}\E[\1_{\{T_{1}\leq\sigma\}}g(x_{i})\Pi_{1}^{i}|\Pi_0=\pi]}.
\label{beben}
\end{eqnarray}
Since $R_{0}$ is deterministic and by using Lemma \ref{lemme-def-op-H}, we recognize that the first term of the right hand side of equation (\ref{beben}) is $Hg(\pi,R_{0})$.
We now turn to the second term of the right hand side of equation (\ref{beben}) which is given by
\begin{eqnarray*}
\E[\1_{\{S_{1}\leq R_{0}\}}\sum_{i=1}^{q}g(x_{i})\Pi_{1}^{i}|\Pi_0=\pi] & = &
\E[v_{N}(\Pi_{1})\1_{\{S_{1}\leq R_{0}\}}|\Pi_0=\pi] \\
&=& \ Gv_{N}(\pi,R_{0}),
\end{eqnarray*}
from Lemma \ref{lemme-tech-event} and the definition of $G$.
Recall that from Definition \ref{def-G-H-I-J-L} one has $J(v_{N},g)=Hg + Gv_{N}$ thus, one obtains
\begin{eqnarray*}
\E[g(X_{\sigma\wedge T_{1}})|\Pi_0=\pi] \ =\ J(v_{N},g)(\pi,R_{0})
&\leq& \sup_{u\geq 0}J(v_{N},g)(\pi,u)\\
&=& L(v_{N},g)(\pi)= v_{N-1}(\pi).
\end{eqnarray*}
Set now $2\leq n\leq N$ and assume that $\E[g(X_{\tau})|\Pi_0=\pi]\leq v_{N-(n-1)}(\pi)$, for all $\tau\in \Sigma_{n-1}^{Y}$.
Proposition \ref{lemme value fonction} yields
\begin{eqnarray*}
 \lefteqn{\E[g(X_{\sigma\wedge T_{n}})|\Pi_0=\pi]}\\
 &=&\sum_{k=0}^{n-1}\sum_{i=1}^{q}\E[\1_{\{ T_{k}\leq \sigma\}}\1_{\{R_{k}<\te_{i}\}}g\circ\Phi(x_{i},R_{k})e^{-\Lambda(x_{i},R_{k})}\Pi_{k}^{i}|\Pi_0=\pi]\\
&&+\sum_{i=1}^{q}\E[\1_{\{T_{n}\leq \sigma\}}g(x_{i})\Pi_{n}^{i}|\Pi_0=\pi].
\end{eqnarray*}
As in the case $n=1$, the term for $k=0$ equals $Hg(\pi,R_{0})$. Notice that 
for $k\geq 1$, $\1_{\{T_{k}\leq \sigma\}}=\1_{\{T_{k}\leq \sigma\}} \1_{\{T_{1}\leq \sigma\}}$ and that
$\1_{\{T_{1}\leq \sigma\}}=\1_{\{S_{1}\leq R_{0}\}}$ is $\Fo_{T_{1}}$-measurable.
By taking the conditional expectation w.r.t. $\Fo_{T_{1}}$ it follows that $\E[\Xi\1_{\{S_{1}\leq R_{0}\}}|\Pi_0=\pi]=\E[\Xi |\Pi_0=\pi]$
where $\Xi$ is defined by
\begin{eqnarray*} 
\Xi&=&\E\Big[\sum_{k=1}^{n-1}\sum_{i=1}^{q}\1_{\{T_{k}\leq \sigma\}}\1_{\{R_{k}<\te_{i}\}}g\circ\Phi(x_{i},R_{k})e^{-\Lambda(x_{i},R_{k})}\Pi_{k}^{i}\\
&&\qquad+\sum_{i=1}^{q}\1_{\{ T_{n}\leq \sigma\}}g(x_{i})\Pi_{n}^{i}\big|\Fo_{T_{1}}\Big].
\end{eqnarray*}
Therefore,  we obtain
\begin{equation} \label{eq-theo-dyn1}
\E[g(X_{\sigma\wedge T_{n}})|\Pi_0=\pi]=
 Hg(\pi,R_{0})+\E[\Xi\1_{\{S_{1}\leq R_{0}\}}|\Pi_0=\pi].
\end{equation}
We now use the Markov property of the chain $(\Pi_{k})_{k\geq 0}$. Indeed, for $k\geq 1$, one has $\Pi_{k}=\Pi_{k-1}\circ \theta$, where $\theta$ is the translation operator of the $(\Fo_{T_{n}})_{n\in\N}$-Markov chain $(\Pi_{n},Y_{n},S_{n})_{n\in \N}$.
Moreover, when $T_{1}\leq \sigma$, one has, from Proposition \ref{propB4}, $R_{k}=\widetilde{R}^{1}_{k-1}\circ\theta$ (indeed, we pointed out in Remark \ref{rq-Rn-Rnbar} that $R_{k}$ can be replaced by $\overline{R}_{k}$ defined in Lemma \ref{lemmeB2}) and $\sigma=T_{1}+\widetilde{\sigma}\circ\theta$ where $\widetilde{R}^{1}_{k-1}$ and $\widetilde{\sigma}$ are defined in Definition \ref{defB3} and Proposition \ref{propB4} (with $l=1$ in the present case).
Since for $k\geq 1$, $T_{k}=T_{1}+T_{k-1}\circ \theta$, one has
$\1_{\{T_{k}\leq \sigma\}}=\1_{\{T_{k-1}\leq \widetilde{\sigma}\}}\circ\theta$.
Finally, combining the Markov property of the chain $(\Pi_{k})_{k\geq 0}$ and Proposition \ref{lemme value fonction}
we have $\Xi=w(\Pi_{1})$ with 
$ w(\pi)=\E[g(X_{\widetilde{\sigma}\wedge T_{n-1}})|\Pi_0=\pi]$.
Moreover, one has $w(\pi)\leq v_{N-(n-1)}(\pi)$ from the induction assumption since $\widetilde{\sigma}\wedge T_{n-1}\in\Sigma_{n-1}^{Y}$ (indeed, both $\widetilde{\sigma}$ and $T_{n-1}$ are $(\Fo_{t})_{t\geq 0}$-stopping times from Corollary \ref{corB6} and Lemma \ref{prop-Tn-stop-time} respectively). One has then 
\begin{equation}\label{eq-theo-dyn2}
\Xi\leq v_{N-(n-1)}(\Pi_{1}).
\end{equation}
Finally, combining Eq.~\eqref{eq-theo-dyn1} and~\eqref{eq-theo-dyn2}, one has
$$\E[g(X_{\sigma\wedge T_{n}})|\Pi_0=\pi] \leq  Hg(\pi,R_{0})+\E[v_{N-(n-1)}(\Pi_{1})\1_{\{S_{1}\leq R_{0}\}}|\Pi_0=\pi].$$
In the second term, we recognize the operator $G$ and one has
\begin{eqnarray*}
\E[g(X_{\sigma\wedge T_{n}})|\Pi_0=\pi]
&\leq&  Hg(\pi,R_{0})+Gv_{N-(n-1)}(\pi,R_{0})\\
&=& J(v_{N-(n-1)},g)(\pi,R_{0})\\
&\leq& \sup_{u\geq 0} J(v_{N-(n-1)},g)(\pi,u)\\
&=& L(v_{N-(n-1)},g)(\pi) \ =\  v_{N-n}(\pi),
\end{eqnarray*}
that proves the induction.
\findemo\\

We now prove the reverse inequality by constructing a sequence of $\epsilon$-optimal stopping times.
\begin{definition}\label{def-Sn-epsilon}For $\epsilon>0$, $1\leq n\leq N$ and for $\pi\in\PR(E_{0})$, we define
$$r_{n}^{\epsilon}(\pi)=\inf\left\{u>0 \text{ : } J(v_{N-n},g)(\pi,u)>v_{N-n-1}(\pi)-\epsilon\right\}.$$
Consider $R_{1,0}^{\epsilon}=r_{0}^{\epsilon}(\Pi_{0})$ and for $2\leq n\leq N$,
$$\left\{\begin{array}{ll}
R_{n,0}^{\epsilon}&=r_{n-1}^{\epsilon/ 2}(\Pi_{0}),\\
R_{n,k}^{\epsilon}&=r_{n-1-k}^{\epsilon/(2^{k+1})}(\Pi_{k})\1_{\{R_{n,k-1}^{\epsilon}\geq S_{k}\}}\text{ for $1\leq k\leq n-2$},\\
R_{n,n-1}^{\epsilon}&=r_{0}^{\epsilon/(2^{n-1})}(\Pi_{n-1})\1_{\{R_{n,n-2}^{\epsilon}\geq S_{n-1}\}},
\end{array}\right.$$
and finally set
$$U_{n}^{\epsilon}=\sum_{k=1}^{n}R_{n,k-1}^{\epsilon}\wedge S_{k}.$$
\end{definition}
The following lemma describes the effect of the translation operator $\theta$ on the sequence $(R_{n,k}^{\epsilon})_{1\leq n\leq N,0\leq k\leq n-1}$.
\begin{lemme}\label{Rnk-epsilon-theta}
For $n\geq 2$ and $1\leq k\leq n-1$, on the set $\{T_{1}\leq U_{n}^{2\epsilon}\}$, one has 
$$R_{n-1,k-1}^{\epsilon}\circ\theta=R_{n,k}^{2\epsilon}.$$ 
\end{lemme}
\demo
For $n=2$, one just has to prove that on the event $\{T_{1}\leq U_{2}^{2\epsilon}\}$, one has $R_{1,0}^{\epsilon}\circ\theta=R_{2,1}^{2\epsilon}$. Yet, from the definition of the sequence $(R_{n,k}^{\epsilon})_{1\leq n\leq N,0\leq k\leq n-1}$, one has
$R_{1,0}^{\epsilon}\circ\theta= r_{0}^{\epsilon}(\Pi_{1})$ and $R_{2,1}^{2\epsilon}= r_{0}^{\frac{2\epsilon}{2}}(\Pi_{1})\1_{\{R_{2,0}^{2\epsilon}\geq S_{1}\}}$.
The result follows since we are on the event $\{T_{1}\leq U_{2}^{2\epsilon}\}=\{R_{2,0}^{2\epsilon}\geq S_{1}\}$.
For a fixed $n\geq 3$, we prove the lemma by induction on $1\leq k\leq n-1$. Set $k=1$. One has from the definition on the sequence $(R_{n,k}^{\epsilon})_{1\leq n\leq N,0\leq k\leq n-1}$,
$R_{n-1,0}^{\epsilon}\circ \theta= r_{n-2}^{\frac{\epsilon}{2}}(\Pi_{1})$ and $R_{n,1}^{2\epsilon}= r_{n-2}^{\frac{2\epsilon}{4}}(\Pi_{1})\1_{\{R_{n,0}^{2\epsilon}\geq S_{1}\}}$.
We obtain $R_{n-1,0}^{\epsilon}\circ \theta=R_{n,1}^{2\epsilon}$ because we have assumed that we are on the event $\{T_{1}\leq U_{n}^{2\epsilon}\}=\{R_{n,0}^{2\epsilon}\geq S_{1}\}$. The propagation of the induction is similar to the case $k=1$.
\findemo\\

Equipped with this preliminary result, we may now prove that $(U_{n}^{\epsilon})_{1\leq n\leq N}$ is a sequence of $\epsilon$-optimal stopping times with respect to the filtration.
generated by the observations.%
\begin{theorem}\label{theo-Sn-epsilon} 
For all $1\leq n\leq N$ and $\epsilon>0$, one has $U_{n}^{\epsilon}\in \Sigma^{Y}_{n}$ and 
$$\E[g(X_{U_{n}^{\epsilon}})|\Pi_0=\pi]\geq v_{N-n}(\pi)-\epsilon.$$
\end{theorem}
\demo Let $n\in \{1,\ldots,N\}$. First notice that, as a direct consequence of Proposition \ref{propB5}, $U_{n}^{\epsilon}$ is an $(\Fo_{t})_{t\geq 0}$-stopping time since, by construction, the $R_{n,k}^{\epsilon}$ are $\Fo_{T_{k}}$-measurable and satisfy the condition $R_{n,k}^{\epsilon}=0$ on the event $\{S_{k}>R_{n,k-1}^{\epsilon}\}$. It is also clear that $U_{n}^{\epsilon}\leq\sum_{k=1}^{n} S_{k} = T_{n}$. Thus, one has $U_{n}^{\epsilon}\in \Sigma^{Y}_{n}$.
Let us now prove the second assessment by induction. Set $n=1$. Let $\pi\in\PR(E_{0})$, we denote $r_{0}^{\epsilon}=r_{0}^{\epsilon}(\pi)$. Since $R_{1,0}^{\epsilon}=r_{0}^{\epsilon}$ is deterministic, one has clearly $R_{1,0}^{\epsilon}\in \Sigma^{Y}$.
Consequently, by using the same arguments as in the proof of Proposition \ref{theo-value-fonction}, we obtain
\begin{align*}
\E[g(X_{R_{1,0}^{\epsilon}\wedge S_{1}})|\Pi_0=\pi]=&
Hg(\pi,r_{0}^{\epsilon})+Gv_{N}(\pi,r_{0}^{\epsilon})=J(v_{N},g)(\pi,r_{0}^{\epsilon}).
\end{align*}
Finally, the definition of $r_{0}^{\epsilon}$ yields $J(v_{N},g)(\pi,r_{0}^{\epsilon})\geq v_{N-1}(\pi)-\epsilon$ thus one has
$$\E[g(X_{R_{1,0}^{\epsilon}\wedge S_{1}})|\Pi_0=\pi]\geq v_{N-1}(\pi)-\epsilon.$$
Now set $2\leq n\leq N$ and assume that 
$\E[g(X_{U_{n-1}^{\epsilon}})|\Pi_0=\pi]\geq v_{N-(n-1)}(\pi)-\epsilon$, for all $\epsilon>0$.
Proposition \ref{lemme value fonction} yields
\begin{eqnarray*}
\lefteqn{ \E[g(X_{U_{n}^{2\epsilon}})|\Pi_0=\pi]}\\
&=&\sum_{k=0}^{n-1}\sum_{i=1}^{q}\E\left[\1_{\{T_{k}\leq U_{n}^{2\epsilon}\}}\1_{\{R_{n,k}^{2\epsilon}<\te_{i}\}}g\circ\Phi(x_{i},R_{n,k}^{2\epsilon})e^{-\Lambda(x_{i},R_{n,k}^{2\epsilon})}\Pi_{k}^{i}\big|\Pi_0=\pi\right]\\
&&+\sum_{i=1}^{q}\E[\1_{\{T_{n}\leq U_{n}^{2\epsilon}\}}g(x_{i})\Pi_{n}^{i}|\Pi_0=\pi].
\end{eqnarray*}
Denote $r_{n-1}^{\epsilon}=r_{n-1}^{\epsilon}(\pi)$. As in the case $n=1$, the term for $k=0$ equals $Hg(\pi,r_{n-1}^{\epsilon})$ since $R_{n,0}^{2\epsilon}=r_{n-1}^{\epsilon}(\Pi_{0})$. Take the conditional expectation w.r.t. $\Fo_{T_{1}}$ in the other terms. One has then,
\begin{equation}\label{eq-theo-dyn-eps3}
\E[g(X_{U_{n}^{2\epsilon}})|\Pi_0=\pi]=
 Hg(\pi,r_{n-1}^{\epsilon})+\E[\Xi'\1_{\{T_{1}\leq U_{n}^{2\epsilon}\}}|\Pi_0=\pi],
\end{equation}
with
\begin{eqnarray*}
 \Xi'&=&\E\Big[\sum_{k=1}^{n-1}\sum_{i=1}^{q}\1_{\{T_{k}\leq U_{n}^{2\epsilon}\}}\1_{\{R_{n,k}^{2\epsilon}<\te_{i}\}}g\circ\Phi(x_{i},R_{n,k}^{2\epsilon})e^{-\Lambda(x_{i},R_{n,k}^{2\epsilon})}\Pi_{k}^{i}\\
&&\qquad+\sum_{i=1}^{q}\1_{\{ T_{n}\leq U_{n}^{2\epsilon}\}}g(x_{i})\Pi_{n}^{i}\big|\Fo_{T_{1}}\Big].
\end{eqnarray*}
Our objective is to apply the Markov property of $(\Pi_{k})_{k\in\N}$ in the term $\Xi'$. Recall that, from Lemma \ref{Rnk-epsilon-theta}, one has $R_{n-1,k-1}^{\epsilon}\circ\theta=R_{n,k}^{2\epsilon}$ for $n\geq 2$ and $1\leq k\leq n-1$ on the event $\{T_{1}\leq U_{n}^{2\epsilon}\}=\{S_{1}\leq R_{n,0}^{2\epsilon}\}$ (the equality of these events stems from Lemma~\ref{lemme-tech-event}). Thus, on this set one has
\begin{align*}
U_{n}^{2\epsilon}&=S_{1}+\sum_{k=2}^{n}R_{n,k-1}^{2\epsilon}\wedge S_{k} = T_{1}+\sum_{k=2}^{n}(R_{n-1,k-2}^{\epsilon}\circ\theta)\wedge (S_{k-1}\circ\theta)\\
&=T_{1}+U_{n-1}^{\epsilon}\circ\theta.
\end{align*}
Besides, recall that $T_{k}=T_{1}+T_{k-1}\circ\theta$, for $k\geq1$.
Consequently, on the set $\{T_{1}\leq U_{n}^{2\epsilon}\}$, one has
$\1_{\{T_{k}\leq U_{n}^{2\epsilon}\}}=\1_{\{T_{k-1}\leq U_{n-1}^{\epsilon}\}}\circ\theta$
and thus, combining the Markov property of the chain $(\Pi_{k})_{k\geq 0}$ and Proposition \ref{lemme value fonction},
we have
$$\Xi'(\Pi_1)=w'(\Pi_1),$$ with 
%
$w'(\pi)=\E\left[g(X_{U_{n-1}^{\epsilon}})\big|\Pi_0=\pi\right]$.
Moreover, thanks to the induction assumption, one has $w'(\pi)\geq v_{N-(n-1)}(\pi)-\epsilon$ so that one obtains
\begin{equation}\label{eq-theo-dyn-eps4}
\Xi'\geq v_{N-(n-1)}(\Pi_{1})-\epsilon.
\end{equation}
Finally, combining equation~\eqref{eq-theo-dyn-eps3} and~\eqref{eq-theo-dyn-eps4} and noticing that, according to Lemma~\ref{lemme-tech-event}, $\{T_{1}\leq U_{n}^{2\epsilon}\}=\{S_{1}\leq r_{n-1}^{\epsilon}\}$, one obtains
\begin{eqnarray*} 
{\E[g(X_{U_{n}^{2\epsilon}})|\Pi_0=\pi]}
&\geq& Hg(\pi,r_{n-1}^{\epsilon})+\E[v_{N-(n-1)}(\Pi_{1})\1_{\{S_{1}\leq r_{n-1}^{\epsilon}\}}|\Pi_0=\pi]-\epsilon\\
&=& J(v_{N-(n-1)},g)(\pi,r_{n-1}^{\epsilon})-\epsilon\\
&\geq& v_{N-n}(\pi)-2\epsilon,
\end{eqnarray*}
from the definition of $r_{n-1}^{\epsilon}$, showing the result.
\findemo
%
\section{Numerical approximation by quantization}\label{section-quantif}
In this section, we are interested in the computational issue for our optimal stopping problem under partial observation. Indeed, we want to compute a numerical approximation of the value function (\ref{opt-stop-pb}) and propose a computable $\epsilon$-optimal stopping time.

As we have seen in the previous section, the value function $v$ can be obtained by iterating the dynamic programming operator $L$.
However, the operator $L$ involves conditional expectations that are in essence difficult to compute and iterate numerically. We manage to overcome this difficulty by combining two special properties of our problem. On the one hand,  the underlying process $(\Pi_n,S_n)$ in the expression of the operator $L$ is a Markov chain. Therefore, it can be discretized using
a quantization technique which is a powerful method suitable for numerical computation and iteration of conditional expectations. On the other hand, the recursion on the functions
$(v_n)_{0\leq n\leq N}$ involving the operator $L$ can be transformed into a recursion on suitably defined random variables. Thus they are easier to iterate numerically as we do not need to compute an approximation of each $v_n$ on the whole state space.

This section is organized as follows. We first explain how the recursion on the functions $(v_n)_{0\leq n\leq N}$ can be transformed into a recurrence on random variables involving only the Markov chain $(\Pi_n,S_n)$. Then, we present a quantization technique to discretize this Markov chain. Afterwards, we construct a discretized version of the main operators in Definition~\ref{def op chapeau} that is used to build an approximation of the value function in Definition~\ref{def v chap}, and a computable $\epsilon$-optimal stopping time. The main results of this section are Theorems~\ref{theo-conv} and  \ref{th arret chap} that prove the convergence of our approximation scheme and provide a rate of convergence.


\bigskip

We first explain how the dynamic programming equations on the functions $(v_n)_{0\leq n\leq N}$ yield a recursion on the random variables $\big(v_{n}(\Pi_{n})\big)_{0\leq n\leq N}$.
Introduce now the sequence $(V_{n})_{0\leq n\leq N}$ of random variables defined by
$$V_{n}=v_{n}(\Pi_{n}).$$
In other words, one has 
\begin{eqnarray}
V_{N}&=&\sum_{i=1}^{q}g(x_{i})\Pi_{N}^{i},\label{def Vn}\\
V_{n}&=&\sup_{u\geq 0} \E\big[\sum_{i=1}^qg\circ \Phi(x_i,u)\Pi^i_n\1_{\{u<t^*_i\}}\1_{\{S_{n+1}>u\}}+V_{n+1}\1_{\{S_{n+1}\leq u\}}|\Pi_{n}\big],\nonumber
\end{eqnarray}
for $0\leq n\leq N-1$.
Notice that $V_N$ is known and the expression of $V_n$ involves only $V_{n+1}$ and the Markov chain $(\Pi_n,S_n)$. Thus, the sequence $(V_n)_{0\leq n\leq N}$ is completely characterized by the system~(\ref{def Vn}). In addition, $V_0=v_0(\Pi_0)=v(\Pi_0)$. Thus to approximate the value function $v$ at the initial point of our process, it is sufficient to provide an approximation of the sequence of random variables $(V_N)_{0\leq n\leq N}$.

%
\subsection{The quantization approach}
There exists an extensive literature on quantization methods for random variables and processes. We do not pretend to present here an exhaustive panorama of these methods. However, the interested reader may for instance, consult the following works \cite{bally03,gray98,pages04} and references therein.
Consider $X$ an $\mathbb{R}^r$-valued random variable such that $\| X \|_p < \infty$ where $\| X \|_p$ denotes the $L^{p}$-nom of $X$: $\| X \|_p=( \mathbb{E}[|X|^{p}])^{1/p}$. 
Let $\nu$ be a fixed integer, the optimal $L^{p}$-quantization of the random variable $X$ consists in finding the best possible $L^{p}$-approximation of $X$ by a random vector $\widehat{X}$ taking at most $\nu$ values: $\widehat{X}\in \{x^{1},\ldots,x^{\nu}\}$.
This procedure consists in the following two steps:
\begin{enumerate}
\item Find a finite weighted grid $\Gamma\subset \mathbb{R}^{r}$ with $\Gamma= \{x^{1},\ldots,x^{\nu}\}$.
\item Set $\widehat{X}=\widehat{X}^{\Gamma}$ where $\widehat{X}^{\Gamma}=proj_{\Gamma}(X)$ with $proj_{\Gamma}$ denotes the closest neighbour projection on $\Gamma$.
\end{enumerate}
The asymptotic properties of the $L^{p}$-quantization are given by the following result, see e.g. \cite{pages04}.

\begin{theorem}
\label{theore}
If $\mathbb{E}[|X|^{p+\eta}]<+\infty$ for some $\eta>0$ then one has
\begin{eqnarray*}
\lim_{\nu\rightarrow \infty} \nu^{p/r} \min_{|\Gamma|\leq \nu} \| X-\widehat{X}^{\Gamma}\|^{p}_{p}& =& J_{p,r} 
\left(\int |h|^{r/(r+p)}(u)du\right)^{1+p/r},
\end{eqnarray*}
where the distribution of $X$ is $P_{X}(du)=h(u) \lambda_{r}(du)+\mu$ with $\mu\perp\lambda_{r}$, $J_{p,r}$ a constant and $ \lambda_{r}$ the Lebesgue measure in $\mathbb{R}^{r}$.
\end{theorem}
There exists a similar procedure for the optimal quantization of a Markov chain. 
Our approximation method is based on the quantization of the Markov chain $(\Pi_{k},S_{k})_{k\leq N}$. Thus, from now on, we will denote, for $0\leq k\leq N$,
$\Theta_k=(\Pi_{k},S_{k})$.
The CLVQ (Competitive Learning Vector Quantization) algorithm \cite[Section 3]{bally03} provides for each time step $0\leq k\leq N$ a finite grid $\Gamma_k$ of $\PR(E_{0})\times\R^+$ as well as the transition matrices $(\widehat{Q}_k)_{0\leq k\leq N-1}$ from $\Gamma_k$ to $\Gamma_{k+1}$. Let~$p~\geq~1$ such that for all $k\leq N$, $\Pi_k$ and $S_k$ have finite moments at least up to order~$p$ and let $proj_{\Gamma_{k}}$ be the nearest-neighbor projection from $\PR(E_{0})\times\R^+$ onto $\Gamma_k$. The quantized process $(\widehat{\Theta}_k)_{k\leq N}=(\widehat{\Pi}_{k},\widehat{S}_{k})_{k\leq N}$ with value for each $k$ in the finite grid $\Gamma_k$ of $\PR(E_{0})\times\R^+$ is then defined by
\begin{equation*}
(\widehat{\Pi}_{k},\widehat{S}_{k})=proj_{\Gamma_{k}}(\Pi_{k},S_{k}).\label{Z_n hat Z_n-mesurable}
\end{equation*}
We will also denote by $\Gamma_k^{\Pi}$, the projection of $\Gamma_k$ on $\PR(E_{0})$, and by $\Gamma_k^{S}$, the projection of $\Gamma_k$ on $\R^{+}$.

Some important remarks must be made concerning the quantization. On the one hand, the optimal quantization has nice convergence properties stated by Theorem \ref{theore}. Indeed, the $L^p$-quantization error $\|\Theta_k-\widehat{\Theta}_k\|_p$ goes to zero when the number of points in the grids goes to infinity. However, on the other hand, the Markov property is not maintained by the algorithm and the quantized process is generally not Markovian. Although the quantized process can be easily transformed into a Markov chain
, this chain will not be homogeneous. It must be pointed out that the quantized process $(\widehat{\Theta}_k)_{k\in\N}$ depends on the starting point $\Theta_0$ of the process.

In practice, we begin with the computation of the quantization grids, 
which merely requires to be able to simulate the process. Notice that in our case, what is actually simulated is the sequence of observation $(Y_{k},S_{k})_{0\leq k\leq N}$. We are then able to compute the filter $(\Pi_{k})_{0\leq k\leq N}$ thanks to the recursive equation provided by Proposition \ref{prop-rec-filtre}. The grids are only computed once and for all and may be stored off-line. Our schemes are then based on the following simple idea: we replace the process by its quantized approximation within the different recursions. The computation is thus carried out in a very simple way since the quantized process has finite state space.
%
\subsection{Approximation of the value function}
Our approximation scheme of the sequence $(V_{n})_{0\leq n\leq N}$ follows the same lines as in \cite{saporta10}, but once more, the results therein cannot be applied directly as the Markov chain $(\Theta_k)_{k\in\N}$ is not the underlying Markov chain of some PDMP. Our approach decomposes in two steps. The first one will be to discretize the time-continuous maximization of the operator $L$ to obtain a maximization over a finite set. The second step consists in replacing the Markov chain $(\Theta_{n})_{n\in\N}=(\Pi_{n},S_{n})_{n\in\N}$ by its quantized approximation $(\widehat \Theta_{n})_{n\in\N}=(\widehat \Pi_{n},\widehat S_{n})_{n\in\N}$ within the dynamic programming equation. Thus, the conditional expectations will become easily tractable finite sums.

Let us first build a finite time grid to discretize the continuous-time maximization in the expression of the operator $L$. The maximum is originally taken over the set $[0,\infty[$. However, it can be seen from Definition~\ref{def-G-H-I-J-L} that $J(v,h)(\pi,u)=J(v,h)(\pi,t^*_q)$ for all $u\geq t^*_q$. Indeed, the random variable $S_1$ is bounded by the greatest deterministic exit time $t^{*}_q$ that is finite thanks to Assumption~\ref{hyp-ts-bounded}. Therefore, the maximization set can be reduced to the compact set $[0,t^*_q]$. Instead of directly discretizing the set $[0,t^*_q]$, we will actually discretize the subsets $]t^*_m,t^*_{m+1}[$. The reason why we want to exclude the points $t^*_m$ from our grid is technical and will be explained with Lemma~\ref{Xi3-indicator}. Now, it seems natural to distinguish wether $\ts_{m}=\ts_{m+1}$ or $\ts_{m}<\ts_{m+1}$.
\begin{definition} 
Let $M\subset \{0,\ldots,q-1\}$ be the set of indices $m$ such that $\ts_{m}<\ts_{m+1}$.
\end{definition}
Notice that $M$ is not empty because it contains at least the index 0 since we assumed that $\ts_{1}>0=\ts_0$. We can now build our approximation grid.

%
%
\begin{definition} Let $\Delta>0$ be such that
\begin{equation}\label{eq-condition-Delta}
\Delta<\frac{1}{2}\min\left\{|\te_{i}-\te_{j}|\text{ with } 0\leq i,j,\leq q \text{ such that }\te_{i}\neq \te_{j}\right\}.
\end{equation}
For all $m\in M$, let $Gr_{m}(\Delta)$ be the finite grid on $]\ts_{m};\ts_{m+1}[$ defined as follows 
$$Gr_{m}(\Delta)=\{\ts_{m}+i\Delta,1\leq i\leq i_{m}\}\union \{\ts_{m+1}-\Delta\},$$ 
where $i_{m}=\max \{i\in\N \text{ such that } \ts_{m}+i\Delta\leq \ts_{m+1}-\Delta \}$. We also denote
$Gr(\Delta)=\union_{ m\in M}Gr_{m}(\Delta).$
\end{definition}
\begin{remarque}\label{rq-Grm}Let $m\in M$. Notice that, thanks to Eq.~\eqref{eq-condition-Delta}, $Gr_{m}(\Delta)$ is not empty. Moreover, it satisfies two properties that will be crucial in the sequel:
\begin{description}
\item[a.]{for all $t\in [\ts_{m};\ts_{m+1}]$, there exists $u\in Gr_{m}(\Delta)$ such that $|u-t|\leq \Delta$,}
\item[b.]{for all $u\in Gr_{m}(\Delta)$ and $0<\eta<\Delta$, one has $[u-\eta;u+\eta]\subset ]\ts_{m};\ts_{m+1}[$.}
\end{description}
\end{remarque}
A discretized maximization operator $L^{d}$ is then defined as follows.
\begin{definition}\label{def K}
Let $L^{d}$: $B(\PR(E_{0}))\times B(E) \rightarrow B(\PR(E_{0}))$ be defined for all $\pi\in\PR(E_{0})$ by
$$L^{d}(v,h)(\pi) = \max_{ m\in M}\big\{\max_{u\in Gr_{m}(\Delta)}\{ J(v,h)(\pi,u)\}\big\}\vee Kv(\pi),$$
with $Kv(\pi)= J(v,h)(\pi,t^*_q)=Gv(\pi,\ts_{q})=\E[v(\Pi_{1})| \Pi_0=\pi].$
\end{definition}
We now proceed to our second step: replacing the Markov chain $(\Theta_{n})_{n\in\N}=(\Pi_{n},S_{n})_{n\in\N}$ by its quantized approximation $(\widehat \Theta_{n})_{n\in\N}=(\widehat \Pi_{n},\widehat S_{n})_{n\in\N}$ within the operators involved in the construction of the value function.
\begin{definition}\label{def op chapeau}
We define the \emph{quantized operators} $\widehat{G}_{n}$, $\widehat{H}_{n}$, $\widehat{J}_{n}$,  $\widehat{K}_{n}$ and $\widehat L^{d}_{n}$ for $n\in\{1,\ldots,N\}$, $v\in B(\Gamma_{n})$, $h\in B(E)$, $\pi\in \Gamma_{n-1}^{\Pi}$ and $u\geq 0$ as follows
\begin{align*}
\widehat{G}_{n}v(\pi,u)=&\hspace{0.3cm}\E[v(\widehat{\Pi}_{n})\1_{\{\widehat{S}_{n}\leq u\}}| \widehat \Pi_{n-1}=\pi],\\
\widehat{H}_{n}h(\pi,u)=&\hspace{0.3cm}\sum_{i=1}^{q}\pi^{i}\1_{\{u< \ts_{i}\}}h\circ\Phi(x_{i},u)\E[\1_{\{\widehat{S}_{n}> u\}}| \widehat \Pi_{n-1}=\pi],\\
\widehat{J}_{n}(v,h)(\pi,u)=&\hspace{0.3cm}\widehat{H}_{n}h(\pi,u)+\widehat{G}_{n}v(\pi,u),\\
\widehat{K}_{n}v(\pi)=&\hspace{0.3cm}\widehat{J}_{n}(v,h)(\pi,t^*_q)=\hspace{0.3cm}\E[v(\widehat{\Pi}_{n})| \widehat \Pi_{n-1}=\pi],\\
\widehat L^{d}_{n}(v,h)(\pi)=&\hspace{0.3cm}\max_{m\in M}\big\{\max_{u\in Gr_{m}(\Delta)}\{ \widehat J_{n}(v,h)(\pi,u)\}\big\}\vee \widehat K_{n}v(\pi).
\end{align*}
\end{definition}
The quantized approximation of the value functions naturally follows.
\begin{definition}\label{def v chap}
For $0\leq n\leq N$, define the functions $\widehat v_{n}$ on $\Gamma_{n}^{\Pi}$ as follows
$$\left\{\begin{array}{rll}
\widehat{v}_{N}(\pi)&=\sum_{i=1}^{q}g(x_{i})\pi^{i}&\text{for all $\pi\in \Gamma_{N}^{\Pi}$,}\\
\widehat{v}_{n-1}(\pi)&=\widehat{L}_{n}^{d}(\widehat{v}_{n},g)(\pi)&\text{for all $\pi\in \Gamma_{n-1}^{\Pi}$ and $1\leq n\leq N$.}
\end{array}\right.$$
For $0\leq n\leq N$, let
$\widehat V_{n}=\widehat v_{n}(\widehat \Pi_{n}).$
\end{definition}
We may now state our main result for the numerical approximation.
\begin{theorem}\label{theo-conv} Suppose that for all $0\leq n\leq N-1$,
\begin{equation}\label{conditionDelta}
\Delta>(2C_{\lambda})^{-1/2}\|S_{n+1}-\widehat S_{n+1}\|_{p}^{1/2},
\end{equation}
then, one has the following bound for the approximation error
\begin{eqnarray*}
\|V_{n}-\widehat V_{n}\|_{p} &\leq& \|V_{n+1}-\widehat V_{n+1}\|_{p}+a\Delta+ b\|S_{n+1}-\widehat S_{n+1}\|_{p}^{\frac{1}{2}}\\
&&+ c_{n}\|\Pi_{n}-\widehat \Pi_{n}\|_{p}+2[v_{n+1}]\|\Pi_{n+1}-\widehat \Pi_{n+1}\|_{p},
\end{eqnarray*}
where $a=[g]_{2}+2C_{g}C_{\lambda}$, $b=2
C_{g}(2C_{\lambda})^{\frac{1}{2}}$ and $c_{n}=[v_{n}]+4C_{g}+2[v_{n+1}]$ with $[v_{n}]$, $[v_{n+1}]$ defined in Proposition \ref{prop-lip-v} and $[g]_{2}$ defined in Assumption \ref{hyp-g-lip}.
\end{theorem}
Theorem \ref{theo-conv} establishes the convergence of our approximation scheme and provides a bound for the rate of convergence. More precisely, it gives a rate for the $L^{p}$ convergence of $\widehat V_{0}$ towards $V_{0}$. Indeed, one has $\|V_{N}-\widehat V_{N}\|_{p}=\|\sum_{i=1}^{q}g(x_{i})\big(\Pi_{N}^{i}-\widehat \Pi_{N}^{i}\big)\|_{p}\leq C_{g}\|\Pi_{N}-\widehat \Pi_{N}\|_{p}$, so by virtue of Theorem \ref{theo-conv} $\big|V_{0}-\widehat V_{0}\big|$ can be made arbitrarily small when the quantization errors $(\|\Theta_{n}-\widehat \Theta_{n}\|_{p})_{0\leq n\leq  N}$ go to zero i.e. when the number of points in the quantization grids goes to infinity.

In order to prove Theorem \ref{theo-conv}, we proceed similarly to \cite{saporta10} and split the approximation error into four terms 
$\|V_{n}-\widehat V_{n}\|_{p}\leq \Xi_{1}+\Xi_{2}+\Xi_{3}+\Xi_{4}$,
with
\begin{align*}
\Xi_{1}&=\|v_{n}(\Pi_{n})-v_{n}(\widehat \Pi_{n})\|_{p},\\
\Xi_{2}&=\|L(v_{n+1},g)(\widehat \Pi_{n})-L^{d}(v_{n+1},g)(\widehat \Pi_{n})\|_{p},\\
\Xi_{3}&=\|L^{d}(v_{n+1},g)(\widehat \Pi_{n})-\widehat L^{d}_{n+1}(v_{n+1},g)(\widehat \Pi_{n})\|_{p},\\
\Xi_{4}&=\|\widehat L^{d}_{n+1}(v_{n+1},g)(\widehat \Pi_{n})-\widehat L^{d}_{n+1}(\widehat v_{n+1},g)(\widehat \Pi_{n})\|_{p}.
\end{align*}
To obtain bounds for each of these terms, one needs to study the regularity of the operators and the value functions $v_n$. The results are detailed in Appendix~\ref{apx lip}. In particular, we establish in Proposition \ref{prop-lip-v} that the value functions $v_n$ are Lipschitz continuous, yielding a bound for the first term.
%
\begin{lemme}
The first term $\Xi_{1}$ is bounded as follows
$$\|v_{n}(\Pi_{n})-v_{n}(\widehat \Pi_{n})\|_{p}\leq [v_{n}]\|\Pi_{n}-\widehat \Pi_{n}\|_{p}.$$
\end{lemme}
The other error terms are studied separately in the following sections.
%
\subsubsection{Second term of the error}
For the second error term, we investigate the consequences of replacing the continuous maximization in operator $L$ by a discrete one on $Gr(\Delta)$.
\begin{lemme}
For all $m\in M$, $v\in B(\PR(E_{0}))$ and $\pi\in\PR(E_{0})$ one has
$$
\big|\sup_{u\in [\ts_{m};\ts_{m+1}[} J(v,g)(\pi,u)-\max_{u\in Gr_{m}(\Delta)} J(v,g)(\pi,u)\big|
\leq \left([g]_{2}+C_{g}C_{\lambda}+C_{v}C_{\lambda}\right)\Delta.
$$
\end{lemme}
\demo 
We use Definition~\ref{def-op-J-L} to split operator $J$ into a sum of continuous operators $J^m$. Thus, one has  
$$\sup_{u\in [\ts_{m};\ts_{m+1}[} J(v,g)(\pi,u)=\sup_{u\in [\ts_{m};\ts_{m+1}]} J^m(v,g)(\pi,u).$$
The function $u\rightarrow J^{m}(v,h)(\pi,u)$ being continuous, there exists $\overline{t}\in [\ts_{m};\ts_{m+1}]$ such that $\sup_{u\in [\ts_{m};\ts_{m+1}]} J^{m}(v,h)(\pi,u)=J^{m}(v,h)(\pi,\overline{t})$. Moreover, from Remark \ref{rq-Grm}.a, one may chose $\overline{u}\in Gr_{m}(\Delta)$ so that $|\overline{u}- \overline{t}|\leq \Delta$. Propositions \ref{prop-lip-H} and \ref{prop-lip-G} stating the Lipschitz continuity of $J^{m}$ then yield
\begin{align*}
0 &\leq  \sup_{u\in [\ts_{m};\ts_{m+1}]} J^{m}(v,h)(\pi,u)-\max_{u\in Gr_{m}(\Delta)} J^{m}(v,h)(\pi,u)\\
&\leq J^{m}(v,h)(\pi,\overline{t})-J^{m}(v,h)(\pi,\overline{u})\\
&\leq \left([g]_{2}+C_{g}C_{\lambda}+C_{v}C_{\lambda}\right)|\overline{t}-\overline{u}|\leq \left([g]_{2}+C_{g}C_{\lambda}+C_{v}C_{\lambda}\right)\Delta,
\end{align*}
showing the result.
\findemo
\begin{lemme}The second term $\Xi_{2}$ is bounded as follows
$$\|L(v_{n+1},g)(\widehat \Pi_{n})-L^{d}(v_{n+1},g)(\widehat \Pi_{n})\|_{p}\leq \left([g]_{2}+2C_{g}C_{\lambda}\right)\Delta.$$
\end{lemme}
\demo 
This is a straightforward consequence of the previous lemma once it has been noticed that for all $a$, $b$, $c$, $d\in\R$, one has $|a\vee b - c\vee d|\leq |a-c|\vee|b-d|$. Notice also that Proposition \ref{prop-lip-v} provides $C_{v_{n+1}}\leq C_{g}$.
\findemo
%
\subsubsection{Third term of the error}
To investigate the third error term, we use the properties of quantization to bound the error made by replacing an operator by its quantized approximation. As in \cite{saporta10}, we must first deal with non-continuous indicator functions. 
The fact that the $t^*_m$ and a small neighborhood around them do not belong to the discretization grid $Gr(\Delta)$ is crucial to obtain the following lemma.
\begin{lemme}\label{Xi3-indicator}
For all $0\leq n\leq N-1$, $m\in M$ and $0<\eta<\Delta$, one has
$$\big\|\max_{u\in Gr_{m}(\Delta)}\E[|\1_{\{S_{n+1}\leq u\}}-\1_{\{\widehat S_{n+1}\leq u\}}|\widehat \Pi_{n}]\big\|_{p}
\leq \eta^{-1}{\|S_{n+1}-\widehat S_{n+1}\|_{p}}+2\eta C_{\lambda}.$$
\end{lemme}
\demo 
Let $0<\eta<\Delta$. The difference of the indicator functions equals 1 if and only if $S_{n+1}$ and $\widehat S_{n+1}$ are on different  sides of $u$. Therefore, if the difference of the indicator functions equals 1, either $|S_{n+1}-u|\leq \eta$, or $|S_{n+1}-u|> \eta$ and in the latter case $|S_{n+1}-\widehat S_{n+1}|> \eta$ too since $|S_{n+1}-\widehat S_{n+1}|> |S_{n+1}-u|$. One has
$|\1_{\{S_{n+1}\leq u\}}-\1_{\{\widehat S_{n+1}\leq u\}}|\leq \1_{\{|S_{n+1}-\widehat S_{n+1}|> \eta\}}+\1_{\{|S_{n+1}-u|\leq \eta\}},$
leading to
\begin{multline*}
\big\|\max_{u\in Gr_{m}(\Delta)}\E\big[|\1_{\{S_{n+1}\leq u\}}-\1_{\{\widehat S_{n+1}\leq u\}}|\big|\widehat \Pi_{n}\big]\big\|_{p}\\
\leq \|\1_{\{|S_{n+1}-\widehat S_{n+1}|> \eta\}}\|_{p}
+ \big\|\max_{u\in Gr_{m}(\Delta)}\E[ \1_{\{|S_{n+1}-u|\leq \eta\}} |\widehat \Pi_{n}]\big\|_{p}.
\end{multline*}
On the one hand, Markov inequality yields
$$\|\1_{\{|S_{n+1}-\widehat S_{n+1}|> \eta\}}\|_{p}= \PP(|S_{n+1}-\widehat S_{n+1}|> \eta)^{\frac{1}{p}}\leq{\|S_{n+1}-\widehat S_{n+1}\|_{p}}{\eta}^{-1}.$$
On the other hand, since $u\in Gr_{m}(\Delta)$, one has 
$[u-\eta;u+\eta]\subset ]\ts_{m};\ts_{m+1}[$ from Remark \ref{rq-Grm}.b, thus $S_{n+1}$ has an absolutely continuous distribution on the interval $[u-\eta;u+\eta]$ since it does not contain any of the $\ts_{i}$. Besides, recall that $\widehat \Theta_{n}=proj_{\Gamma_{n}}(\Theta_{n})$, hence,
the following inclusions of $\sigma$-fields
$\sigma(\widehat \Pi_{n})\subset\sigma(\widehat \Theta_{n})\subset \sigma(\Theta_{n})$. We also have $\sigma(\Theta_{n})\subset \Fo_{T_{n}}\subset \F_{T_{n}}$, the law of iterated conditional expectations provides
\begin{eqnarray*}
\E[ \1_{\{|S_{n+1}-u|\leq \eta\}} |\widehat \Pi_{n}]
&=&\E\big[\E\big[\E[ \1_{\{|S_{n+1}-u|\leq \eta\}} | \F_{T_{n}}]\big| \Fo_{T_{n}}\big]\big|\widehat \Pi_{n}\big]\\
& \leq &\E\Big[\E[ \int_{u-\eta}^{u+\eta}\lambda\big(\Phi(Z_{n},s)\big)ds\Big| \Fo_{T_{n}}]\Big|\widehat \Pi_{n}\Big]\\
&=&\E\Big[ \sum_{i=1}^{q} \Pi_{n}^{i}\int_{u-\eta}^{u+\eta}\lambda\big(\Phi(x_{i},s)\big)ds\big|\widehat \Pi_{n}\Big].
\end{eqnarray*}
Finally, one obtains
$\E[ \1_{\{|S_{n+1}-u|\leq \eta\}}|\widehat \Pi_{n}]\leq 2\eta C_{\lambda},$
showing the result.
\findemo
\begin{lemme}\label{Xi3-K}
For all $0\leq n\leq N-1$, one has
\begin{multline*}
|Kv_{n+1}(\widehat \Pi_{n})-\widehat K_{n+1}v_{n+1}(\widehat \Pi_{n})|\\
\leq [v_{n+1}]\E\big[|\Pi_{n+1}-\widehat \Pi_{n+1}|\big|\widehat \Pi_{n}\big]+(2C_{g}+2[v_{n+1}])\E\left[|\Pi_{n}-\widehat \Pi_{n}|\big|\widehat \Pi_{n}\right].
\end{multline*}
\end{lemme}
\demo By the definitions of operators $K$ and $\widehat{K}_{n+1}$, one has
\begin{eqnarray}
\lefteqn{|Kv_{n+1}(\widehat \Pi_{n})-\widehat K_{n+1}v_{n+1}(\widehat \Pi_{n})|}\nonumber\\
&=&|\E[v_{n+1}(\Pi_{n+1})| \Pi_{n}=\widehat \Pi_{n} ]-\E[v_{n+1}(\widehat \Pi_{n+1})| \widehat \Pi_{n} ]|\nonumber\\
&\leq&{|\E[v_{n+1}(\Pi_{n+1})| \Pi_{n}=\widehat \Pi_{n} ]-\E[v_{n+1}(\Pi_{n+1})| \widehat \Pi_{n} ]|}\nonumber\\
&&+{|\E[v_{n+1}(\Pi_{n+1})-v_{n+1}(\widehat \Pi_{n+1})| \widehat \Pi_{n} ]|}.\label{eq dif K}
\end{eqnarray}
The second term in the right-hand side of Eq. (\ref{eq dif K}) is readily bounded by using Proposition \ref{prop-lip-v} stating that $v_{n+1}$ is Lipschitz continuous
\begin{equation*}
|\E[v_{n+1}(\Pi_{n+1})-v_{n+1}(\widehat \Pi_{n+1})| \widehat \Pi_{n} ]|\leq[v_{n+1}]\E\big[|\Pi_{n+1}-\widehat \Pi_{n+1}|\big|\widehat \Pi_{n}\big].
\end{equation*} 
To deal with the first term in the right-hand side of Eq. (\ref{eq dif K}), we need to use the special properties of quantization. Indeed,
one has $(\widehat \Pi_{n},\widehat S_{n})=proj_{\Gamma_{n}}(\Pi_{n},S_{n})$ so that 
we have the inclusion of $\sigma$-fields
$\sigma(\widehat \Pi_{n})\subset \sigma(\Pi_{n},S_{n})$.  The law of iterated conditional expectations gives
$$\E[v_{n+1}(\Pi_{n+1})\big| \widehat \Pi_{n} ]=\E\big[\E[v_{n+1}(\Pi_{n+1})|(\Pi_{n},S_{n})]\big| \widehat \Pi_{n} \big].$$
Moreover, Proposition \ref{prop-markov} yields 
$\E[v_{n+1}(\Pi_{n+1})|{(}\Pi_{n},S_{n}{)}]=\E[v_{n+1}(\Pi_{n+1})|\Pi_{n}]$, as the conditional distribution of $\Pi_{n+1}$ w.r.t. $(\Pi_{n},S_{n})$ merely depends on $\Pi_{n}$.
In addition, $|\E[v_{n+1}(\Pi_{n+1})| \Pi_{n}=\widehat \Pi_{n} ]$ is $\sigma(\widehat \Pi_{n})$-measurable.
One has then
\begin{eqnarray*}
\lefteqn{|\E[v_{n+1}(\Pi_{n+1})| \Pi_{n}=\widehat \Pi_{n} ]-\E[v_{n+1}(\Pi_{n+1})| \widehat \Pi_{n} ]|}\\
&=& \left|\E\big[\E[v_{n+1}(\Pi_{n+1})| \Pi_{n}=\widehat \Pi_{n} ]-\E[v_{n+1}(\Pi_{n+1})| \Pi_{n} ]\big|\widehat \Pi_{n}\big]\right|\\
&=&|\E[Kv_{n+1}(\widehat \Pi_{n})-Kv_{n+1}(\Pi_{n})|\widehat \Pi_{n}]|,
\end{eqnarray*}
by definition of $K$. Finally, one has
\begin{eqnarray*}
{\lefteqn{|\E[v_{n+1}(\Pi_{n+1})| \Pi_{n}=\widehat \Pi_{n} ]-\E[v_{n+1}(\Pi_{n+1})| \widehat \Pi_{n} ]|}}\\
&\leq&2(C_g+[v_{n+1}])\E\left[|\Pi_{n}-\widehat \Pi_{n}|\big|\widehat \Pi_{n}\right],
\end{eqnarray*}
thanks to Propositions \ref{prop-lip-K} and \ref{prop-lip-v} stating the Lipschitz continuity of operator $K$ and function $v_{n+1}$.
\findemo
\begin{lemme}\label{Lemme-Xi3}
If $\Delta$ satisfies Condition (\ref{conditionDelta}), a upper bound for the third term $\Xi_{3}$ is
\begin{eqnarray*}
\lefteqn{\|L^{d}(v_{n+1},g)(\widehat \Pi_{n})-\widehat L^{d}_{n+1}(v_{n+1},g)(\widehat \Pi_{n})\|_{p}}\\
&\leq& [v_{n+1}]\|\Pi_{n+1}-\widehat \Pi_{n+1}\|_{p}+(4C_{g}+2[v_{n+1}])\|\Pi_{n}-\widehat \Pi_{n}\|_{p}\\
&&+2C_{g} (2 C_{\lambda})^{1/2}{\|S_{n+1}-\widehat S_{n+1}\|_{p}}^{1/2}.
\end{eqnarray*}
\end{lemme}
\demo 
One has
\begin{eqnarray*}
\lefteqn{|L^{d}(v_{n+1},g)(\widehat \Pi_{n})-\widehat L^{d}_{n+1}(v_{n+1},g)(\widehat \Pi_{n})|}\\
&\leq&\max_{m\in M}\big\{\max_{u\in Gr_{m}(\Delta)} |J(v_{n+1},g)(\widehat \Pi_{n},u)-\widehat J_{n+1}(v_{n+1},g)(\widehat \Pi_{n},u)|\big\}\\
&&\vee |Kv_{n+1}(\widehat \Pi_{n})-\widehat K_{n+1}v_{n+1}(\widehat \Pi_{n})|.
\end{eqnarray*}
The term involving operator $K$ was studied in the previous lemma.
Let us now study the term involving operator $J$. Set $m$ in $M$,
 $u$ in $Gr_{m}(\Delta)$ and define $\alpha(\pi,\pi',s')=\sum_{i=1}^{q}\pi^{i}g\big(\Phi(x_{i},u)\big)\1_{\{s'>u\}}+v_{n+1}(\pi')\1_{\{s'\leq u\}}$. One has then
\begin{eqnarray*}
\lefteqn{|J(v_{n+1},g)(\widehat \Pi_{n},u)-\widehat J_{n+1}(v_{n+1},g)(\widehat \Pi_{n},u)|}\\
&=&|J^{m}(v_{n+1},g)(\widehat \Pi_{n},u)-\widehat J_{n+1}(v_{n+1},g)(\widehat \Pi_{n},u)|\\
&=&\left|\E[\alpha(\Pi_{n},\Pi_{n+1},S_{n+1})|\Pi_{n}=\widehat \Pi_{n}]-\E[\alpha(\widehat \Pi_{n},\widehat \Pi_{n+1},\widehat S_{n+1})|\widehat \Pi_{n}]\right|\leq A+B,
\end{eqnarray*}
where
\begin{align*}
A=&\left|\E[\alpha(\Pi_{n},\Pi_{n+1},S_{n+1})|\Pi_{n}=\widehat \Pi_{n}]-\E[\alpha(\Pi_{n},\Pi_{n+1},S_{n+1})|\widehat\Pi_{n}]\right|,\\
B=&|\E[\alpha(\Pi_{n},\Pi_{n+1},S_{n+1})-\alpha(\widehat \Pi_{n},\widehat \Pi_{n+1},\widehat S_{n+1})\big|\widehat \Pi_{n}]|.
\end{align*}
Using the boundedness of $g$ and $v_{n+1}$ as well as the Lipschitz continuity of $v_{n+1}$ given in Proposition \ref{prop-lip-v}, we get a upper bound for the second term
\begin{eqnarray}\label{eq-bound-A}
B&\leq &C_{g}\E\big[|\Pi_{n}-\widehat \Pi_{n}|\big|\widehat \Pi_{n}\big] + [v_{n+1}]\E\big[|\Pi_{n+1}-\widehat \Pi_{n+1}|\big|\widehat \Pi_{n}\big]\nonumber\\
&&+ 2C_{g}\E\left[|\1_{\{S_{n+1}\leq u\}}-\1_{\{\widehat S_{n+1}\leq u\}}|\big|\widehat \Pi_{n}\right].
\end{eqnarray}
For the first term, we use the properties of quantization as in the previous proof to obtain
\begin{eqnarray*}
A&=&\left|\E\big[\E[\alpha(\Pi_{n},\Pi_{n+1},S_{n+1})|\Pi_{n}=\widehat \Pi_{n}]-\E[\alpha(\Pi_{n},\Pi_{n+1},S_{n+1})|\Pi_{n}]\big| \widehat \Pi_{n}\big]\right|.
\end{eqnarray*}
We now
recognize operator $J^{m}$,
and from Propositions \ref{prop-lip-H} and \ref{prop-lip-G}, one has
\begin{eqnarray}\label{eq-bound-B}
A&=&\E[J^{m}(v_{n+1},g)(\widehat \Pi_{n},u)-J^{m}(v_{n+1},g)(\Pi_{n},u)|\widehat \Pi_{n}]\nonumber\\
&\leq& (3C_{g}+2[v_{n+1}])\E\left[|\widehat{\Pi}_{n}-\Pi_{n}\big|\widehat \Pi_{n}\right].
\end{eqnarray}
We gather the bounds provided by Eq.~\eqref{eq-bound-A} and \eqref{eq-bound-B}  to obtain
\begin{eqnarray}\label{eq J-Jchap}
\lefteqn{|J(v_{n+1},g)(\widehat \Pi_{n},u)-\widehat J_{n+1}(v_{n+1},g)(\widehat \Pi_{n},u)|}\nonumber\\
&\leq&(4C_{g}+2[v_{n+1}])\E\big[|\Pi_{n}-\widehat \Pi_{n}|\big|\widehat \Pi_{n}\big] + [v_{n+1}]\E\big[|\Pi_{n+1}-\widehat \Pi_{n+1}|\big|\widehat \Pi_{n}\big]\nonumber\\
&&+ 2C_{g}\E\left[|\1_{\{S_{n+1}\leq u\}}-\1_{\{\widehat S_{n+1}\leq u\}}|\big|\widehat \Pi_{n}\right].
\end{eqnarray}
Finally, combining the result for operators $J$ and Lemma~\ref{Xi3-K}, we obtain 
\begin{eqnarray*}
\lefteqn{|L^{d}(v_{n+1},g)(\widehat \Pi_{n})-\widehat L^{d}_{n+1}(v_{n+1},g)(\widehat \Pi_{n})|}\\
&\leq&[v_{n+1}]\E\left[|\Pi_{n+1}-\widehat \Pi_{n+1}|\big|\widehat \Pi_{n}\right]+(4C_{g}+2[v_{n+1}])\E\left[|\Pi_{n}-\widehat \Pi_{n}|\big|\widehat \Pi_{n}\right]\\
&&+2C_{g}\max_{u\in Gr(\Delta)}\E\left[|\1_{\{S_{n+1}\leq u\}}-\1_{\{\widehat S_{n+1}\leq u\}}|\big|\widehat \Pi_{n}\right].
\end{eqnarray*}
We conclude by taking the $L^{p}$ norm in the equation above and using Lemma \ref{Xi3-indicator} to bound the last term
\begin{eqnarray*}
\lefteqn{\|L^{d}(v_{n+1},g)(\widehat \Pi_{n})-\widehat L^{d}_{n+1}(v_{n+1},g)(\widehat \Pi_{n})\|_p}\\
&\leq&[v_{n+1}]\|\Pi_{n+1}-\widehat \Pi_{n+1}\|_p+(4C_{g}+2[v_{n+1}])\|\Pi_{n}-\widehat \Pi_{n}\|_p\\
&&+2C_{g}(\eta^{-1}{\|S_{n+1}-\widehat S_{n+1}\|_{p}}+2\eta C_{\lambda}),
\end{eqnarray*}
for some $0<\eta<\Delta$. The best choice for $\eta$ minimizing the error is when $\eta$ satisfies
\begin{equation*}
\eta^{-1}{\|S_{n+1}-\widehat S_{n+1}\|_{p}}=2\eta C_{\lambda},
\end{equation*}
which yields $\eta=(2 C_{\lambda})^{-1/2}({\|S_{n+1}-\widehat S_{n+1}\|_{p}}^{1/2}$. If $\Delta$ satisfies Condition (\ref{conditionDelta}), one has $\eta<\Delta$ as required for this optimal choice.
\findemo
%
\subsubsection{Fourth term of the error}
Finally, the fourth error term is bounded using Lipschitz properties.
\begin{lemme}\label{Lemme-Xi4}
The fourth term $\Xi_{4}$ is bounded as follows
\begin{eqnarray*}
\lefteqn{\|\widehat L^{d}_{n+1}(v_{n+1},g)(\widehat \Pi_{n})-\widehat L^{d}_{n+1}(\widehat v_{n+1},g)(\widehat \Pi_{n})\|_{p}}\\
&\leq& [v_{n+1}]\|\Pi_{n+1}-\widehat \Pi_{n+1}\|_{p}+\|V_{n+1}-\widehat V_{n+1}\|_{p}.
\end{eqnarray*}
\end{lemme}
\demo 
One has
\begin{eqnarray}
\lefteqn{\|\widehat L^{d}_{n+1}(v_{n+1},g)(\widehat \Pi_{n})-\widehat L^{d}_{n+1}(\widehat v_{n+1},g)(\widehat \Pi_{n})\|_{p}}\nonumber\\
&=&\big\|\max_{m\in M}\max_{u\in Gr_{m}(\Delta)}\big\{\widehat{H}_{n+1}g(\widehat \Pi_{n},u)+\widehat{G}_{n+1}v_{n+1}(\widehat \Pi_{n},u)\big\}\vee\widehat K_{n+1}v_{n+1}(\widehat \Pi_{n})\nonumber\\
&&-\max_{m\in M}\max_{u\in Gr_{m}(\Delta)}\big\{\widehat{H}_{n+1}g(\widehat \Pi_{n},u)+\widehat{G}_{n+1}\widehat v_{n+1}(\widehat \Pi_{n},u)\big\}\vee\widehat K_{n+1}\widehat v_{n+1}(\widehat \Pi_{n})\big\|_{p},\nonumber\\
&\leq&\big\|\max_{m\in M}\max_{u\in Gr_{m}(\Delta)}\E\left[\big(v_{n+1}(\widehat{\Pi}_{n+1})-\widehat v_{n+1}(\widehat{\Pi}_{n+1})\big)\1_{\{\widehat{S}_{n+1}\leq u\}}\big| \widehat{\Pi}_{n}\right]\nonumber\\
&&\vee \E[v_{n+1}(\widehat{\Pi}_{n+1})-\widehat v_{n+1}(\widehat{\Pi}_{n+1})| \widehat{\Pi}_{n}] \big\|_{p}\nonumber\\
&\leq&\|v_{n+1}(\widehat{\Pi}_{n+1})-\widehat v_{n+1}(\widehat{\Pi}_{n+1})\|_{p}.\label{Jvchap}
\end{eqnarray}
We now introduce $v_{n+1}(\Pi_{n+1})$ to split this term into two differences. The Lipschitz continuity of $v_{n+1}$ stated by Proposition \ref{prop-lip-v} allows us to bound the first term while we recognize $V_{n+1}$ and $\widehat V_{n+1}$ in the second one.
\begin{eqnarray*}
\lefteqn{\|\widehat L^{d}_{n+1}(v_{n+1},g)(\widehat \Pi_{n})-\widehat L^{d}_{n+1}(\widehat v_{n+1},g)(\widehat \Pi_{n})\|_{p}}\\
&\leq&\|v_{n+1}(\widehat{\Pi}_{n+1})-v_{n+1}(\Pi_{n+1})\|_{p}+\|v_{n+1}(\Pi_{n+1})-\widehat v_{n+1}(\widehat{\Pi}_{n+1})\|_{p}\\
&\leq& [v_{n+1}]\left\|\Pi_{n+1}-\widehat \Pi_{n+1}\right\|_{p}+\|V_{n+1}-\widehat V_{n+1}\|_{p}.
\end{eqnarray*}
Hence, the result.
\findemo

\subsection{Numerical construction of an $\epsilon$-optimal stopping time}\label{section-numeric-epsilon-opt}
As in the previous section, we follow the idea of \cite{saporta10} and we use both the Markov chain $(\Theta_{n})_{0\leq n\leq N}$ and its quantized approximation $(\widehat{\Theta}_{n})_{0\leq n\leq N}$ to approximate the expression of the $\epsilon$-optimal stopping time introduced in Definition \ref{def-Sn-epsilon}. 
We check that we thus obtain actual stopping times for the observed filtration $(\Fo_{t})_{t\geq 0}$ and that the expected reward when stopping then is a good approximation of the value function $V_{0}$.
For all $(\pi,s) \in \PR(E_{0})\times\R^{+}$ and $0\leq n\leq N$, we denote $(\widehat{\pi}_{n},\widehat{s}_{n})=proj_{\Gamma_{n}}(\pi,s)$. Let 
\begin{equation*}
\widehat{s}^{*}_{N-n}(\pi,s)=\min\{t\in Gr(\Delta):
\widehat{J}_{n}(\widehat{v}_{n},g)(\widehat{\pi}_{n-1},t) = \max_{u\in Gr(\Delta)} \widehat{J}_{n}(\widehat{v}_{n},g)(\widehat{\pi}_{n-1},u)\}.
\end{equation*}
For $1\leq n\leq N$ and $\pi\in\PR(E_{0})$, we define
$$\widehat{r}_{N-n}(\pi,s)=
\left\{\begin{array}{ll}
\ts_{q} & \text{ if } \widehat{K}_{n}\widehat{v}_{n}(\widehat{\pi}_{n-1})> \max_{u\in Gr(\Delta)} \widehat{J}_{n}(\widehat{v}_{n},g)(\widehat{\pi}_{n-1},u),\\
\widehat{s}^{*}_{N-n}(\pi,s) &\text{ otherwise.}
\end{array}\right.
$$
Let now for $n\geq 1$,
$$\left\{\begin{array}{ll}
\widehat{R}_{n,0}&=\widehat{r}_{n-1}(\Pi_{0},S_{0}),\\
\widehat{R}_{n,k}&=\widehat{r}_{n-1-k}(\Pi_{k},S_{k})\1_{\{\widehat{R}_{n,k-1}\geq S_{k}\}}\text{ for $1\leq k\leq n-2$,}
\end{array}\right.$$
and set
$\widehat{U}_{n}=\sum_{k=1}^{n}\widehat{R}_{n,k-1}\wedge S_{k}.$
The following result is a direct consequence of Proposition \ref{propB5}. It is a very strong result as it states that the numerically computable random variables $\widehat{U}_{n}$ are actual $(\Fo_{t})_{t\geq 0}$-stopping times.
\begin{theorem}
For $0\leq n\leq N$, $\widehat{U}_{n}$ is an $(\Fo_{t})_{t\geq 0}$-stopping time. 
\end{theorem}
We now intend to prove that stopping at time $\widehat{U}_{N}$ provides a good approximation of the value function $V_{0}$. For all $\pi \in\PR(E_{0})$ and $0\leq n\leq N$ we therefore introduce the performance when abiding by the stopping rule~$(\widehat{U}_{n})_{0\leq n\leq N}$ and the corresponding random variables
$$\overline{v}_{n}(\pi)=\E[g(X_{\widehat{U}_{N-n}})|\Pi_{0}=\pi],\qquad \overline{V}_{n}=\overline{v}_{n}(\Pi_{n}).$$
\begin{theorem}\label{th arret chap} 
Suppose that for all $0\leq n\leq N-1$,
$$\Delta>({2C_{\lambda}})^{-1/2}{\|S_{n+1}-\widehat S_{n+1}\|_{p}}^{1/2},$$
one has then the following bound for the error between the expected reward when stopping at time $\widehat{U}_{n}$ and the value function
\begin{eqnarray*}
\|V_{n}-\overline V_{n}\|_{p} &\leq& \left\|V_{n+1}-\overline V_{n+1}\right\|_{p}+\left\|V_{n}-\widehat V_{n}\right\|_{p}+\left\|V_{n+1}-\widehat V_{n+1}\right\|_{p}\\
&&+ d_{n}\|\Pi_{n}-\widehat \Pi_{n}\|_{p}+2[v_{n+1}]\|\Pi_{n+1}-\widehat \Pi_{n+1}\|_{p}\\
&&+b\|S_{n+1}-\widehat S_{n+1}\|_{p}^{1/2},
\end{eqnarray*}
where $b=2C_{g}\big(2C_{\lambda}\big)^{1/2}$, $d_{n}=7C_{g}+4[v_{n+1}]$, $[v_{n+1}]$ defined in Proposition~\ref{prop-lip-v}.
\end{theorem}
It is important to notice that $\overline{v}_{N}(\pi)=\sum_{i=1}^{q}g(x_{i})\pi^{i}=v_{N}(\pi)$ and thus $\overline{V}_{N}=V_{N}$. Therefore, the previous theorem proves that $|V_{0}-\overline{V}_{0}|$ goes to zero when the quantization errors $(\|\Theta_{n}-\widehat{\Theta}_{n}\|_{p})_{0\leq n\leq N}$ go to zero. In other words, the expected reward $\overline{V}_{0}$ when stopping at the random time $\widehat{U}_{N}$ can be made arbitrarily close to the value function $V_{0}$ of the partially observed optimal stopping problem \eqref{opt-stop-pb} and hence $\widehat{U}_{N}$ is an $\epsilon$-optimal stopping time.
\demo
The first step consists in finding a recursion satisfied by the sequence $(\overline{V}_{n})_{0\leq n\leq N}$ in order to compare it with the dynamic programming equation giving $(\widehat{V}_{n})_{0\leq n\leq N}$. 
Let $0\leq n\leq N-1$. First of all, Proposition~\ref{lemme value fonction} gives
\begin{eqnarray*} 
\lefteqn{\E[g(X_{\widehat{U}_{N-n}})|\Pi_0]}\\
&=&\sum_{k=0}^{N-n-1}\sum_{i=1}^{q}\E[\1_{\{T_{k}\leq \widehat{U}_{N-n}\}}\1_{\{\widehat{R}_{N-n,k}<\te_{i}\}}g\circ\Phi(x_{i},\widehat{R}_{N-n,k})e^{-\Lambda(x_{i},\widehat{R}_{N-n,k})}\Pi_{k}^{i}|\Pi_0]\\
&&+\sum_{i=1}^{q}\E[\1_{\{T_{N-n}\leq \widehat{U}_{N-n}\}}g(x_{i})\Pi_{n}^{i}|\Pi_0].
\end{eqnarray*}
The term corresponding to $k=0$ in the above sum equals $Hg(\Pi_{0},\widehat{R}_{N-n,0})$. Taking the conditional expectation w.r.t. $\Fo_{T_{1}}$ in the other terms and noticing that one has $\{T_{1}\leq \widehat{U}_{N-n}\}=\{S_{1}\leq \widehat{R}_{N-n,0}\}$ yield
$$
\E[g(X_{\widehat{U}_{N-n}})|\Pi_0]=Hg(\Pi_{0},\widehat{R}_{N-n,0})+\E[\Xi''\1_{\{S_{1}\leq \widehat{R}_{N-n,0}\}}|\Pi_0],
$$
with
\begin{eqnarray*} 
\Xi''&=&\E\Big[\sum_{k=1}^{N-n-1}\sum_{i=1}^{q}\1_{\{T_{k}\leq \widehat{U}_{N-n}\}}\1_{\{\widehat{R}_{N-n,k}<\te_{i}\}}g\circ\Phi(x_{i},\widehat{R}_{N-n,k})e^{-\Lambda(x_{i},\widehat{R}_{N-n,k})}\Pi_{k}^{i}\\
&&\qquad+\sum_{i=1}^{q}\1_{\{ T_{N-n}\leq \widehat{U}_{N-n}\}}g(x_{i})\Pi_{n}^{i}\big|\Fo_{T_{1}}\Big].
\end{eqnarray*}
We now make use of the Markov property of the sequence $(\Pi_{n})_{n\in\N}$ in the term $\Xi''$. Similarly to Lemma \ref{Rnk-epsilon-theta}, for $n\geq 1$, on the set $\{T_{1}\leq \widehat{U}_{N-n}\}$, one has $\widehat{R}_{N-n-1,k-1}\circ\theta=\widehat{R}_{N-n,k}$ for all $1\leq k\leq n-1$. Thus, on the set $\{T_{1}\leq \widehat{U}_{N-n}\}$, one has $\widehat{U}_{N-n}=T_{1}+\widehat{U}_{N-n-1}\circ\theta$. Recall that $\1_{\{T_{k}\leq \widehat{U}_{N-n}\}}=\1_{\{T_{k-1}\leq \widehat{U}_{N-n-1}\}}\circ\theta$.
We may therefore apply the Markov property. Using Proposition \ref{lemme value fonction}, we now obtain $\Xi''=\overline{v}_{n+1}(\Pi_{1})$. Finally, we have 
$$\overline{v}_{n}(\Pi_{0})=Hg(\Pi_{0},\widehat{R}_{N-n,0})+G\overline{v}_{n+1}(\Pi_{0},\widehat{R}_{N-n,0})=J(\overline{v}_{n+1},g)(\Pi_{0},\widehat{R}_{N-n,0}).$$
Recall that $\widehat{R}_{N-n,0}=\widehat{r}_{N-n-1}(\Pi_{0},S_{0})$ and apply the translation operator $\theta^{n}$ to obtain the following recursion
$$\overline{V}_{n}=J(\overline{v}_{n+1},g)(\Pi_{n},\widehat{r}_{N-n-1}(\Pi_{n},S_{n})).$$
We are now able to study the error between $\overline{V}_{n}$ and $\widehat{V}_{n}$. Let us recall that, from its definition, $\widehat{r}_{N-n-1}(\Pi_{n},S_{n})$ equals either $\widehat{s}^{*}_{N-n-1}(\Pi_{n},S_{n})$ or $\ts_{q}$. In the latter case, notice that $J(\overline{v}_{n+1},g)(\Pi_{n},\ts_{q})=K\overline{v}_{n+1}(\Pi_{n})$. Eventually, one has
\begin{eqnarray}\label{prog dyn Vbar}
|\overline{V}_{n}-\widehat{V}_{n}|\leq \1_{\{\widehat{r}_{N-n-1}(\Pi_{n},S_{n})=\ts_{q}\}}A+\1_{\{\widehat{r}_{N-n-1}(\Pi_{n},S_{n})=\widehat{s}^{*}_{N-n-1}(\Pi_{n},S_{n})\}}B
\end{eqnarray}
with
$$\left\{\begin{array}{ll}
A&=|K\overline{v}_{n+1}(\Pi_{n})- \widehat{K}_{n+1}\widehat{v}_{n+1}(\widehat{\Pi}_{n})|,\\
B&=| J(\overline{v}_{n+1},g)(\Pi_{n},\widehat{s}^{*}_{N-n-1}(\Pi_{n}))-\max_{u\in Gr(\Delta)} \widehat{J}_{n+1}(\widehat{v}_{n+1},g)(\widehat{\Pi}_{n},u)|.
\end{array}\right.$$
To bound the first term $A$, we introduce the function $v_{n+1}$. One has
\begin{eqnarray*}
A&\leq &{|K\overline{v}_{n+1}(\Pi_{n})- Kv_{n+1}(\Pi_{n})|} + {|Kv_{n+1}(\Pi_{n})- Kv_{n+1}(\widehat{\Pi}_{n})|}\\
&&+{|Kv_{n+1}(\widehat{\Pi}_{n})- \widehat{K}_{n+1}v_{n+1}(\widehat{\Pi}_{n})|} + {|\widehat{K}_{n+1}v_{n+1}(\widehat{\Pi}_{n})- \widehat{K}_{n+1}\widehat{v}_{n+1}(\widehat{\Pi}_{n})|}\\
&\leq& (a)+(b)+(c)+(d).
\end{eqnarray*}
Let us study these four terms one by one. By definition of $K$, the first term $(a)$
 is bounded by $\E\left[|\overline{V}_{n+1}-V_{n+1}|\big|\Pi_{n}\right]$. For the second term $(b)$, we use Proposition \ref{prop-lip-K} stating the Lipschitz continuity of the operator $K$. The term third term $(c)$ is bounded by Lemma \ref{Xi3-K} and a upper bound of the  fourth term $(d)$ is given by Eq. (\ref{Jvchap}). Thus, one obtains
 \begin{eqnarray*}
 A&\leq&\left\|V_{n+1}-\overline V_{n+1}\right\|_{p}+\left\|V_{n+1}-\widehat V_{n+1}\right\|_{p}
+ 4(C_g+[v_{n+1}])\|\Pi_{n}-\widehat \Pi_{n}\|_{p}\\
&&+2[v_{n+1}]\|\Pi_{n+1}-\widehat \Pi_{n+1}\|_{p}.
 \end{eqnarray*}
We now turn to the second term $B$. In the following computations, denote $\widehat{s}^{*}=\widehat{s}^{*}_{N-n-1}(\Pi_{n},S_{n})$. Its definition yields
$B=| J(\overline{v}_{n+1},g)(\Pi_{n},\widehat{s}^{*})-\widehat{J}_{n+1}(\widehat{v}_{n+1},g)(\widehat{\Pi}_{n},\widehat{s}^{*})|$.
We split this expression into four differences again. On the set $\{\widehat{r}_{N-n-1}(\Pi_{n},S_{n})=\widehat{s}^{*}\}$, one has the equality $J(v_{n+1},g)(\Pi_{n},\widehat{s}^{*})=V_{n+1}$. Hence, one this set, one obtains from Eq. (\ref{prog dyn Vbar})
$$
| J(\overline{v}_{n+1},g)(\Pi_{n},\widehat{s}^{*})-J(v_{n+1},g)(\Pi_{n},\widehat{s}^{*})|
\leq |\overline{V}_{n+1}-V_{n+1}|.
$$
For the other terms, we use Propositions \ref{prop-lip-H} and \ref{prop-lip-G} for the Lipschitz continuity of $J$ and Eq. (\ref{eq J-Jchap}) and (\ref{Jvchap}) to obtain
 \begin{eqnarray*}
B&\leq&\left\|V_{n+1}-\overline V_{n+1}\right\|_{p}+\left\|V_{n+1}-\widehat V_{n+1}\right\|_{p}\\
&&+ (7C_g+4[v_{n+1}])\|\Pi_{n}-\widehat \Pi_{n}\|_{p}+2[v_{n+1}]\|\Pi_{n+1}-\widehat \Pi_{n+1}\|_{p}\\
&&+2C_g(2C_{\lambda})^{1/2}\|S_{n+1}-\widehat S_{n+1}\|_{p}^{1/2},
 \end{eqnarray*}
after optimizing $\eta$. The result is obtained by taking the maximum between $A$ and $B$.
\findemo
%
\section{Numerical example}\label{section-example}
\begin{figure}[thbp]
\begin{center}
\includegraphics[width=8cm]{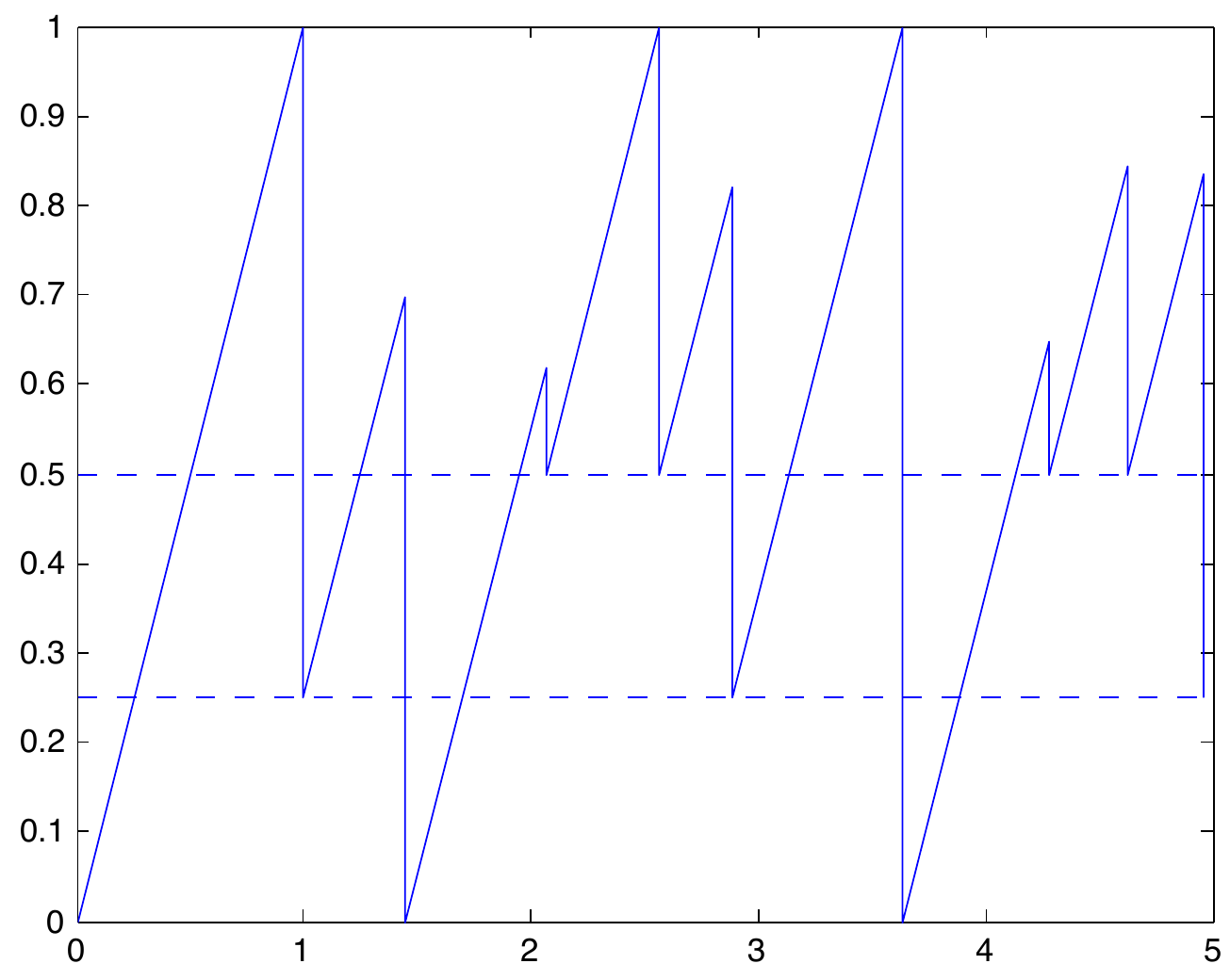}
\caption{A trajectory of the process drawn until the $9^{\text{th}}$ jump time with $a=3$, $v=1$ and $E_{0}=\{0;\frac{1}{4}; \frac{1}{2}\}$. The dotted lines represent the possible post-jump values.}
\label{plotecg}
\end{center}
\end{figure}

We apply our procedure to a simple PDMP similar to the one studied in~\cite{saporta10}.
Let $E=[0;1[$. For $x\in E$ and $t\geq 0$, the flow is defined by $\Phi(x,t)=x+vt$ so that $\ts(x)=(1-x)/v$. We set the jump rate to $\lambda(x)=ax$ for some $a>0$ and the transition kernel $Q(x,\cdot)$ to the uniform distribution on a finite set $E_{0}\subset E$. Thus, the process evolves toward 1 and the closer it gets to 1, the more likely it will jump back to some point of $E_{0}$. A trajectory is represented in Figure \ref{plotecg}.
The observation process is $Y_{n}=\varphi(Z_{n})+W_{n}$ where $\varphi(x)=x$ and $W_{n}\sim \mathcal{N}(0,\sigma^{2})$ for some $\sigma^{2}>0$. Finally, we choose the reward function $g(x)=x$. Our assumptions thus clearly hold.
Simulations are run with $a=3$, $v=1$, $E_{0}=\{0,1/4,1/2\}$, $\sigma^{2}=0.25$ and $N=9$. 
The numerical approximation is implemented as follows. First, we make an exact simulator for the sequence $(Z_n,S_n)$. From the values of $(Z_n)$, one builds the observation sequence $(Y_n)$ that allows for a recursive computation of the filter process thanks to Proposition \ref{prop-rec-filtre}. Thus, we can simulate trajectories of the Markov chain $(\Pi_n,S_n)$ that we feed into the CLVQ algorithm to obtain quantization grids. By Monte Carlo simulations, we can also estimate the quantization errors. To run our numerical procedure, one then needs to choose the parameter $\Delta$ satisfying conditions (\ref{eq-condition-Delta}) and (\ref{conditionDelta}). In this special case, they boil down to
\begin{equation*}
6^{-1/2}\max_{0\leq n\leq N-1}\|S_{n+1}-\widehat S_{n+1}\|_{p}^{1/2}<\Delta<\frac{1}{8}.
\end{equation*}
We have chosen $\Delta$ just above the Monte Carlo approximation of the lower bound. The values are given in the second column of Table \ref{ecg_table} for different grids sizes.

Then, we recursively compute the approximated value functions $\widehat{v}_n$ on the quantization grids. The conditional expectations are now merely weighted sums. The approximation we obtain for the value function of the partially observed optimal stopping problem are given in the fourth column of Table~\ref{ecg_table}.

Finally, we implemented the construction of our $\epsilon$-optimal stopping time and ran $10^6$ Monte Carlo simulations to compute its mean performance. The results are given in the third column of Table~\ref{ecg_table}.

The exact value of $V_{0}$ is unknown but one has as in \cite{saporta10},
\begin{equation}\label{B-empir}
\overline V_{0}=\E[g(X_{\widehat{U}_{N}})] \leq V_{0}=\sup_{\sigma\in\Sigma_{N}^{Y}}\E[g(X_{\sigma})] \leq \E\big[ \sup_{0\leq t \leq T_{N}}g(X_{t})\big].
\end{equation}
Both the first and the last term may be estimated by Monte Carlo simulations. One has thus, with $10^{6}$ trajectories, $\E[ \sup_{0\leq t \leq T_{N}}g(X_{t})]=0.9944$. 
The theoretical bound $B_{th}$ of the error $|V_{0}-\widehat V_{0}|$ provided by Theorem \ref{theo-conv} is computed using the approximated quantization errors. This bound decreases as the number of points in the quantization grids increases, as expected. Moreover, we computed the empirical bound given by Eq.~\eqref{B-empir}
$B_{em}= \max \big\{ |\overline V_{0}-\widehat V_{0}| , |\E[ \sup_{0\leq t \leq T_{N}}g(X_{t})]-\widehat V_{0}| \big\}$.
\begin{table}[tbp]
\begin{center}
\begin{tabular}{cccccc}
\text{Quantization grids} & $\Delta$ & $\overline V_{0}$ &$\widehat V_{0}$& $B_{em}$& $B_{th}$\\
\hline
\text{50 points}  &0.1179&0.7900&0.8135 & 0.181&683\\
\hline
\text{100 points}  &0.0970&0.8031&0.8250& 0.169&467 \\
\hline
\text{300 points}  &0.0731&0.8182&0.8407& 0.154&271 \\
\hline
\text{500 points}  &0.0634&0.8250&0.8477& 0.147&211 \\
\hline
\text{1000 points} &0.0535&0.8313&0.8545& 0.140&152 \\
\hline
\text{2000 points} &0.0453&0.8361&0.8599& 0.135&110 \\
\hline
\text{4000 points} &0.0381&0.8408&0.8643& 0.130&80 \\
\hline
\text{6000 points} &0.0345&0.8430&0.8666& 0.128&67 \\
\hline
\text{8000 points} &0.0321&0.8479&0.8725&0.122&58\\
\hline
\text{10000 points} &0.0303&0.8497&0.8742&0.120&53\\
\hline
\text{12000 points} &0.0290&0.8521&0.8771&0.117&49\\
\end{tabular}
\end{center}
\caption{Simulation results. The terms $B_{em}$ and $B_{th}$ respectively denote an empirical bound and the theoretical bound provided by Theorem \ref{theo-conv} for the error $|V_{0}-\widehat V_{0}|$.}
\label{ecg_table}
\end{table}
%
\appendix
%
\section{Properties of the $(\Fo_{t})_{t\geq 0}$-stopping times} \label{STOP-time}
In this section, we study the special structure of $(\Fo_{t})_{t\geq 0}$-stopping times. 
\begin{lemme}\label{prop-Tn-stop-time}
For all $n\in \N$, $T_{n}$ is an $(\Fo_{t})_{t\geq 0}$-stopping time.
\end{lemme}
\demo Notice that for all $n\in\N$, $\PP(Y_{n}= Y_{n+1})=0$. This stems from the absolute continuity of the distribution of the random variables $(W_{n})_{n\in\N}$ since
$$\left\{Y_{n}=Y_{n+1}\right\}\subset \union_{1\leq i,j \leq q} \left\{W_{n}-W_{n+1}=\varphi(x_{i})-\varphi(x_{j})\right\}.$$
Hence, for all $n\in\N$ and $t\in\R^{+}$, one has $\PP$ a.s. $\left\{T_{n}\leq t\right\}=\left\{N_{t}\geq n\right\}$ where we denote $N_{t}=\sum_{0\leq s\leq t} \1_{\{Y_{s}\neq Y_{s^{-}}\}}$. The process $(N_{t})_{t\geq 0}$ is $\Fo$-adapted thus $\{N_{t}\geq n\}\in \Fo_{t}$ and since the filtration $\Fo$ contains the $\PP$-null sets, one has $\left\{T_{n}\leq t\right\}\in \Fo_{t}$. For all $n\in\N$, $T_{n}$ is therefore an $(\Fo_{t})_{t\geq 0}$-stopping time.
\findemo\\

We now recall Theorem A2 T33 from \cite{bremaud81} concerning the structure of the stopping times for point processes and apply it in our case.
\begin{definition}Define the filtration $(\F_{t}^{p})_{t\geq 0}$ as follows
$$\F_{t}^{p}=\sigma\left(\1_{\{Y_{n}\in A\}}\1_{\{T_{n}\leq s\}}; n\geq 1, 0\leq s\leq t, A\in\mathcal{B}(\R^{d})\right).$$
\end{definition}
\begin{theorem}\label{theo-bremaud}
Let $\sigma$ be an $(\F_{t}^{p})_{t\geq 0}$-stopping time. For all $n\in\N$, there exists a $\F_{T_{n}}^{p}$-measurable non negative random variable $R_{n}$,  such that one has
$$\sigma\wedge T_{n+1} = \big(T_{n}+R_{n}\big)\wedge T_{n+1}\quad on\quad \{\sigma\geq T_{n}\}.$$
\end{theorem}
Our observation process $(Y_{t})_{t\geq 0}$ being a point process that fits the framework developed in \cite{bremaud81}, we apply this Theorem to $(\Fo_{t})_{t\geq 0}$-stopping times.
\begin{proposition}
For all $t\geq 0$, one has $\Fo_{t}=\F_{t}^{p}$.
\end{proposition}
\demo First prove that $\Fo_{t}\subset\F_{t}^{p}$. Let $A\in \mathcal{B}(\R^{d})$ and $0\leq s\leq t$, one has
$$\left\{Y_{s}\in A\right\}=\union_{n\in\N}\big(\left\{T_{n}\leq s< T_{n+1}\right\}\inter\left\{Y_{n}\in A\right\}\big)\in\F_{s}^{p}\subset \F_{t}^{p}.$$
Indeed, in the above equation, we used that $T_{0}$ and $Y_{0}$ are assumed to be deterministic.
For the reverse inclusion, let $A\in \mathcal{B}(\R^{d})$, $n\geq 1$ and $0\leq s\leq t$. Recall that $Y_{n}=Y_{T_{n}}$. One has $\left\{Y_{T_{n}}\in A\right\}\in \Fo_{T_{n}}$ since $(Y_{t})_{t\geq 0}$ is $\Fo$-adapted and $T_{n}$ is an $(\Fo_{t})_{t\geq 0}$-stopping time from Lemma \ref{prop-Tn-stop-time}. Therefore, one has
$\left\{Y_{n}\in A\right\}\inter \left\{T_{n}\leq s\right\}\in\Fo_{s}\subset\Fo_{t}$,
showing the result.
\findemo\\

We may therefore apply Theorem \ref{theo-bremaud} to $(\Fo_{t})_{t\geq 0}$-stopping times. 
\begin{theorem}\label{theo-bremaud-adapted}
Let $\sigma$ be an $(\Fo_{t})_{t\geq 0}$-stopping time. For all $n\in\N$, there exists a non negative random variable $R_{n}$, $\Fo_{T_{n}}$-measurable such that one has
$$\sigma\wedge T_{n+1} = \big(T_{n}+R_{n}\big)\wedge T_{n+1}\quad on\quad \{\sigma\geq T_{n}\}.$$
\end{theorem}
We outline the following result, which is a direct consequence of the above theorem, because it will be used several times in our derivation.
\begin{lemme} \label{lemme-tech-event}
Let $\sigma$ be an $(\Fo_{t})_{t\geq 0}$-stopping time and $(R_{n})_{n\in \N}$ be the sequence of random variables associated to $\sigma$ as introduced in
Theorem \ref{theo-bremaud-adapted}.
For all $n\in \N$,
$\{T_{n}\leq\sigma<T_{n+1}\}=\{T_{n}\leq\sigma\}\inter\{S_{n+1}>R_{n}\}.$
\end{lemme}
\demo Theorem \ref{theo-bremaud-adapted} states that on the event $\{T_{n}\leq\sigma\}$, on has $\sigma\wedge T_{n+1}= T_{n} + (R_{n}\wedge S_{n+1})$ so that, still on the event $\{T_{n}\leq\sigma\}$, one has $(\sigma<T_{n+1}) \Leftrightarrow (R_{n} < S_{n+1})$.
We deduce the result from this observation.
\findemo\\

We now investigate the effect of the translation operator of the Markov chain $(\Pi_{n},Y_{n},S_{n})_{n\in\N}$ on the $(\Fo_{t})_{t\geq 0}$-stopping times. Proposition \ref{prop-markov} states that $(\Pi_{n},Y_{n},S_{n})_{n\in\N}$ is a $(\Fo_{T_{n}})_{n\in\N}$-Markov chain. Let us consider its canonical space $\Omega=(\PR(E_{0})\times \R^{d}\times \R^{+})^{\N}$. Thus, for $\omega=(\omega_{0},\omega_{1},\ldots)\in\Omega$, one has $(\Pi_{n},Y_{n},S_{n})(\omega)=\omega_{n}$. Besides, we define the \emph{translation operator} 
$$\theta : 
\left\{\begin{array}{clc}
\Omega &\rightarrow &\Omega\\
(\omega_{0},\omega_{1},\ldots)&\rightarrow&(\omega_{1},\omega_{2},\ldots)
\end{array}\right. 
$$
We then define $\theta^{0}=Id_{\Omega}$ and recursively for $l\geq 2$, $\theta^{l}=\theta\circ\theta^{l-1}$. Thus, for all $n,l\in\N$, one has $(\Pi_{n},Y_{n},S_{n})\circ \theta^{l}=(\Pi_{n+l},Y_{n+l},S_{n+l})$. As $T_{0}=0$, one has 
$$T_{n}\circ\theta^{l}=\sum_{k=1}^{n}S_{k}\circ\theta^{l}=\sum_{k=1}^{n}S_{k+l}=T_{n+l}-T_{l}.$$
The next results of this section 
are given without proof because their proofs follow the very same lines as in \cite{saporta10} from which they are adapted. However, notice that the results from \cite{saporta10} cannot be applied directly to our case because the sequence $(\Pi_{n},Y_{n},S_{n})_{n\in\N}$, although it is a Markov chain, is not the underlying Markov chain of some PDMP.
Set now $\sigma\in\Sigma^{Y}$. From Theorem \ref{theo-bremaud-adapted}, for all $n\in\N$, there exists a non negative $\Fo_{T_{n}}$-measurable random variable $R_{n}$, such that, on the event $\{\sigma\geq T_{n}\}$, one has
$\sigma\wedge T_{n+1} = \big(T_{n}+R_{n}\big)\wedge T_{n+1}.$
\begin{lemme}\label{lemmeB2}
Let $\sigma$ be an $(\Fo_{t})_{t\geq 0}$-stopping time and $(R_{n})_{n\in \N}$ be the sequence of random variables associated to $\sigma$ as introduced in
Theorem \ref{theo-bremaud-adapted}.
Let $\overline{R}_{0}=R_{0}$ and for $k\geq 1$, $\overline{R}_{k}=R_{k}\1_{\{S_{k}\leq \overline{R}_{k-1}\}}$. One has then
$$\sigma=\sum_{n=1}^{\infty}\overline{R}_{n-1}\wedge S_{n}.$$
\end{lemme}
\begin{remarque}\label{rq-Rn-Rnbar}
This lemma proves that in Theorem \ref{theo-bremaud-adapted}, the sequence $(R_{n})_{n\in\N}$ can be replaced by $(\overline{R}_{n})_{n\in\N}$. Therefore, we can assume, without loss of generality that the sequence $(R_{n})_{n\in\N}$ satisfies the following condition: for all $n\in\N$, $R_{n+1}=0$ on the event $\{S_{n+1}>R_{n}\}$.
\end{remarque}
Since $\Fo_{T_{k}}=\sigma(Y_{j},S_{j}, j\leq k)$ and $R_{k}$ is $\Fo_{T_{k}}$-measurable, there exists a sequence of real-valued measurable functions $(r_{k})_{k\in\N}$ defined on $(\R^{d} \times \R^{+})^{k+1}$ such that
$R_{k}=r_{k}(\mathcal{G}_{k})$, 
where $\mathcal{G}_{k}=(Y_{0},S_{0},\ldots,,Y_{k},S_{k})$. 
\begin{definition}\label{defB3}
Let $\sigma$ be an $(\Fo_{t})_{t\geq 0}$-stopping time and $(r_{n})_{n\in \N}$ be the sequence of functions associated to $\sigma$ as introduced in
Remark \ref{rq-Rn-Rnbar}.
Let $l\geq 1$ and $(\widetilde{R}^{l}_{k})_{k\in\N}$ be a sequence of functions defined on $(\R^{d} \times \R^{+})^{l+1}\times \Omega$ by
$\widetilde{R}^{l}_{0}(\gamma,\omega)=r_{l}(\gamma)$
and for $k\geq 1$,
$\widetilde{R}^{l}_{k}(\gamma,\omega)=r_{l+k}(\gamma,\mathcal{G}_{k-1}(\omega))\1_{\{S_{k}\leq \widetilde{R}^{l}_{k-1}\}}(\gamma,\omega).$
\end{definition}
\begin{proposition}\label{propB4}
Let $\sigma$ be an $(\Fo_{t})_{t\geq 0}$-stopping time and $( \overline{R}_{k})_{k\in\N}$ (respectively,  $(\widetilde{R}^{l}_{k})_{k\in\N}$) be the sequence of functions associated to $\sigma$ as introduced in Lemma \ref{lemmeB2} (respectively, in Definition \ref{defB3}).
Assume that $T_{l}\leq \sigma \leq T_{N}$. For all $k\in\N$, one has then $\widetilde{R}^{l}_{k}(\mathcal{G}_{l},\theta^{l})=\overline{R}_{l+k}$ and
$\sigma=T_{l}+\widetilde{\sigma}(\mathcal{G}_{l},\theta^{l})$,
with $\widetilde{\sigma}$ : $\left(\R^{d} \times \R^{+}\right)^{l+1}\times \Omega\rightarrow \R^{+}$ defined as
$\widetilde{\sigma}(\gamma,\omega)=\sum_{n=1}^{N-l}\widetilde{R}^{l}_{n-1}(\gamma,\omega)\wedge S_{n}(\omega).$
\end{proposition}
\begin{proposition}\label{propB5} Let $(U_{n})_{n\in\N}$ be a sequence of non negative random variables such that for all $n$, $U_{n}$ is $\Fo_{T_{n}}$-measurable and $U_{n+1}=0$ on $\{S_{n+1}>U_{n}\}$. We define
$U=\sum_{n=1}^{\infty}U_{n-1}\wedge S_{n}.$
Then $U$ is an $(\Fo_{t})_{t\geq 0}$-stopping time.
\end{proposition}
\begin{corollaire}\label{corB6}
Let $\sigma$ be an $(\Fo_{t})_{t\geq 0}$-stopping time and  $\widetilde{\sigma}$ be the mapping associated to $\sigma$ introduced in Proposition \ref{propB4}.
For all $\gamma\in(\R^{d} \times \R^{+})^{p+1}$, $\widetilde{\sigma}(\gamma,\cdot)$ is a $(\Fo_{t})_{t\geq 0}$-stopping time.
\end{corollaire}
%
\section{Computation of a conditional expectation}\label{appendix-esp-cond}
The objective of this section is to prove the technical Lemma \ref{lemme-tech-esp-cond} used in the proof of Proposition \ref{lemme value fonction}.
%
\begin{lemme}\label{lemme-tech-esp-cond}
For all $k\in \N$, one has
$\E[\1_{\{S_{k+1}>R_{k}\}}|\F_{T_{k}}]=\1_{\{R_{k}<\te(Z_{k})\}}e^{-\Lambda(Z_{k},R_{k})}.$
\end{lemme}
\demo
First recall some results concerning the random variables $(S_{k})_{k\in\N}$, details may be found in \cite{davis93}.
After a jump of the process to the point $z\in E$, the survival function of the time until the next jump is
$$\phi(t,z)=
\left\{\begin{array}{ll}
1 & \text{if } t\leq0,\\
e^{-\Lambda(z,t)} & \text{if } 0\leq t<t^*(z), \\
0 & \text{if }t\geq t^*(z).
\end{array}\right.$$
Define its generalized inverse
$\psi(u,z)=\inf\{t\geq 0\text{ such that } \phi(t,z)\leq u\}$.
Then, for all $k\in\N$, one has $S_{k+1}=\psi(\Upsilon_k,Z_{k})$,
where $\Upsilon_{k}$ are i.i.d. random variables with uniform distribution on $[0;1]$ independent from $\F_{T_{k}}$.
Thus, one has $\E[\1_{\{S_{k+1}>R_{k}\}}|\F_{T_{k}}]=\E[f(\Upsilon_{k},Z_{k},R_{k})|\F_{T_{k}}]$ where $f(u,z,r)=\1_{\{\psi(u,z)>r\}}$. As $(Z_{k},R_{k})$ is $\F_{T_{k}}$-measurable, $\Upsilon_{k}$ is independent from $\F_{T_{k}}$ and $\E[\1_{\{\psi(\Upsilon_k,z)>r\}}]=\1_{\{r<\te(z)\}}e^{-\Lambda(z,r)}$, \cite[Proposition 11.2]{ouvrard04} yields the result.
\findemo
%
\section{Lipschitz properties}\label{apx lip}
In this section, we derive the Lipschitz properties of our operators in order to obtain them for the value functions~$(v_{n})_{0\leq n\leq N}$. 
Similarly to the proof of Proposition \ref{prop-markov}, we first derive the integral form of operators $G$ and $H$.
%
\begin{lemme}\label{lemme-def-op-H} 
For all $h\in B(E)$, $v\in B(\PR(E_{0}))$ and $(\pi,u)\in \PR(E_{0})\times  \R^{+}$, one has
\begin{eqnarray*}
{Gv(\pi,u)}
&=&Iv(\pi,u) + \sum_{i=1}^{q}\pi^{i}\1_{\{\te_{i}\leq u\}}e^{-\Lambda(x_{i},\te_{i})}\\
&&\times\int_{\R^{d}}v\big(\Psi(\pi,y',\ts_{i})\big)\sum_{j=1}^{q}Q\big(\Phi(x_{i},\te_{i}),x_{j}\big)f_{W}(y'-\varphi(x_{j}))dy',\\
Hh(\pi,u)
&=&\sum_{i=1}^{q}\pi^{i}\1_{\{u< \ts_{i}\}}e^{-\Lambda(x_{i},u)}h\circ\Phi(x_{i},u),
\end{eqnarray*}
where
\begin{eqnarray*}
Iv(\pi,u)
&=& \sum_{i=1}^{q}\pi^{i}\int_{0}^{u\wedge \te_{i}}\Big(\lambda\circ\Phi(x_{i},s')e^{-\Lambda(x_{i},s')} \\
&&\times\int_{\R^{d}}v\big(\Psi(\pi,y',s')\big)\sum_{j=1}^{q}Q\big(\Phi(x_{i},s'),x_{j}\big)f_{W}\big(y'-\varphi(x_{j})\big)dy'\Big)ds'.\\
\end{eqnarray*}
\end{lemme}
Now, notice that the functions $Hh(\pi,\cdot)$ and $Gv(\pi,\cdot)$ are not continuous. However, they are c\`adl\`ag with a finite number of jumps. Therefore, they can be rewritten as sums of continuous functions as follows.
\begin{definition} \label{def-op-J-L}
For all $m\in \{0,\ldots,q-1\}$, we define the operators $G^{m}$: $B(\PR(E_{0})) \rightarrow B(\PR(E_{0})\times \R^{+})$ and $H^{m}$: $B(E) \rightarrow B(\PR(E_{0})\times \R^{+})$ as follows
\begin{itemize}
\item{if $u < \ts_{m}$, $G^{m}v(\pi,u)=Gv(\pi,\ts_{m})$ and $H^{m}h(\pi,u)=Hh(\pi,\ts_{m})$,}
\item{if $u \geq \ts_{m}$,
\begin{eqnarray*}
\hspace{-1cm}G^{m}v(\pi,u)&=&Iv(\pi,u\wedge \ts_{m+1})+ \sum_{i=1}^{m}\pi^{i}e^{-\Lambda(x_{i},\te_{i})}\\ 
&&\times\int_{\R^{d}}v\big(\Psi(\pi,y',\ts_{i})\big)\sum_{j=1}^{q}Q\big(\Phi(x_{i},\te_{i}),x_{j}\big)f_{W}\big(y'-\varphi(x_{j})\big)dy',\\
\hspace{-1cm}H^{m}h(\pi,u)&=&\sum_{i=m+1}^{q}\pi^{i}e^{-\Lambda(x_{i},u\wedge \ts_{m+1})}h\circ\Phi(x_{i},u\wedge \ts_{m+1}).
\end{eqnarray*}
}
\end{itemize}
We also define
$J^{m}(v,h)(\pi,u)=H^{m}h(\pi,u)+G^{m}v(\pi,u).$
\end{definition}

\begin{remarque}\label{rq-continu-Hm}
For all $m\in \{0,\ldots,q-1\}$ and for all $h\in B(E)$, $v\in B(\PR(E_{0}))$ and $(\pi,u)\in \PR(E_{0})\times  \R^{+}$, the functions $u\rightarrow G^{m}v(\pi,u)$, $u\rightarrow H^{m}h(\pi,u)$ and $u\rightarrow J^{m}(v,h)(\pi,u)$ are continuous. Moreover, they are constant on $[0;\ts_{m}]$ and on $[\ts_{m+1};+\infty[$ and one has
\begin{eqnarray*}
Gv(\pi,u)&=&\sum_{m=0}^{q-1}\1_{[\ts_{m},\ts_{m+1}[}(u)G^{m}v(\pi,u),\\
Hh(\pi,u)&=&\sum_{m=0}^{q-1}\1_{[\ts_{m},\ts_{m+1}[}(u)H^{m}h(\pi,u),\\
J(v,h)(\pi,u)&=&\sum_{m=0}^{q-1}\1_{[\ts_{m},\ts_{m+1}[}(u)J^{m}(v,h)(\pi,u).
\end{eqnarray*}
\end{remarque}
%
%
%
We now investigate the Lipschitz properties of our operators.
\begin{proposition}\label{prop-lip-H}
For $m\in M$, $\left((\pi,u), (\tilde \pi,\tilde u)\right)\in (\PR(E_{0}) \times \R^{+})^{2}$, one has
$$|H^{m}g(\pi,u)-H^{m}g(\tilde\pi,\tilde u)|\leq C_{g}|\pi-\tilde\pi|+([g]_{2}+C_{g}C_{\lambda})|u-\tilde u|.$$
\end{proposition}
\demo 
Since the function $u\rightarrow H^{m}h(\pi,u)$ is constant on the intervals $[0;\ts_{m}]$ and $[\ts_{m+1};+\infty[$, we may assume that $u$, $\tilde u\in [\ts_{m};\ts_{m+1}]$ so that one has
$H^{m}g(\pi,u)=\sum_{i=m+1}^{q}\pi^{i}e^{-\Lambda(x_{i},u)}g\circ\Phi(x_{i},u),$
and similarly for $H^{m}g(\tilde \pi,\tilde u)$. Then, on the one hand, one has
\begin{eqnarray*}
{\left|H^{m}g(\pi,u)-H^{m}g(\tilde \pi,u)\right|}
&=&\big|\sum_{i=m+1}^{q}\big(\pi^{i}-\tilde \pi^{i}\big)e^{-\Lambda(x_{i},u)}g\circ\Phi(x_{i},u)\big|\\
&\leq& C_{g}\sum_{i=m+1}^{q}|\pi^{i}-\tilde \pi^{i}|.
\end{eqnarray*}
On the other hand, Lemma A.1 in \cite{saporta10} yields
$$
|e^{-\Lambda(x_{i},u)}g\circ\Phi(x_{i},u)-e^{-\Lambda(x_{i},\tilde u)}g\circ\Phi(x_{i},\tilde u)|\leq ([g]_{2}+C_{g}C_{\lambda})|u-\tilde u|,
$$
showing the result.
\findemo\\

The following technical lemma will be useful to derive the Lipschitz properties of the operator $I$. The first part of its proof is adapted from \cite{pham05}.
\begin{lemme}
For all $\pi$, $\tilde \pi\in \PR(E_{0})$ and $m\in M$, one has
$$\sum_{m=0}^{q-1}\int_{\ts_{m}}^{\ts_{m+1}}\int_{\R^{d}}|\Psi(\pi,y',s')-\Psi(\tilde \pi,y',s')|\overline{\Psi}_{m}(\pi,y',s')dy'ds'\leq 2|\pi-\tilde \pi|.$$
\end{lemme}
\demo
Let $s'\in]\ts_{m};\ts_{m+1}[$ and $y'\in\R^{d}$. In the following computation, we denote $\tau=(\pi,y',s')$ and $\tilde \tau=(\tilde \pi,y',s')$, one has
\begin{eqnarray*}
{|\Psi(\tau)-\Psi(\tilde \tau)|\overline{\Psi}_{m}(\tau)}
&=&\sum_{j=1}^{q}\left|\frac{\Psi^{j}_{m}(\tau)}{\overline{\Psi}_{m}(\tau)}-\frac{\Psi^{j}_{m}(\tilde \tau)}{\overline{\Psi}_{m}(\tilde \tau)}\right|\overline{\Psi}_{m}(\tau)\\
&=&\sum_{j=1}^{q}\left|\frac{\Psi^{j}_{m}(\tau)\overline{\Psi}_{m}(\tilde \tau)-\Psi^{j}_{m}(\tilde \tau)\overline{\Psi}_{m}(\tau)}{\overline{\Psi}_{m}(\tilde \tau)}\right|\\
&\leq&\sum_{j=1}^{q}\left|\Psi^{j}_{m}(\tau)-\Psi^{j}_{m}(\tilde \tau)\right|+\sum_{j=1}^{q}\frac{\Psi^{j}_{m}(\tilde \tau)}{\overline{\Psi}_{m}(\tilde \tau)}\left|\overline{\Psi}_{m}(\tau)-\overline{\Psi}_{m}(\tilde \tau)\right|.
\end{eqnarray*}
Notice that $\sum_{j=1}^{q}\Psi^{j}_{m}(\tilde \tau)=\overline{\Psi}_{m}(\tilde \tau)$ so that the second sum above reduces to $|\overline{\Psi}_{m}(\tau)-\overline{\Psi}_{m}(\tilde \tau)|=\sum_{j=1}^{q}\left|\Psi^{j}_{m}(\tau)-\Psi^{j}_{m}(\tilde \tau)\right|$. Finally, one has
\begin{equation*}
|\Psi(\tau)-\Psi(\tilde \tau)|\overline{\Psi}_{m}(\tau)\leq 2\sum_{j=1}^{q}|\Psi^{j}_{m}(\tau)-\Psi^{j}_{m}(\tilde \tau)|.
\end{equation*}
As $\int_{\R^{d}}f_{W}\big(y'-\varphi(x_{j})\big)dy'=1$ and $\sum_{j=1}^{q}Q\big(\Phi(x_{i},s'),x_{j}\big)=1$, one obtains

\begin{eqnarray*}
\lefteqn{\sum_{m=0}^{q-1}\int_{\ts_{m}}^{\ts_{m+1}}\int_{\R^{d}}\big|\Psi(\pi,y',s')-\Psi(\tilde \pi,y',s')\big|\overline{\Psi}_{m}(\pi,y',s')dy'ds'}\\
&\leq&2\sum_{m=0}^{q-1}\int_{\ts_{m}}^{\ts_{m+1}}\sum_{j=1}^{q}\int_{\R^{d}}\left|\Psi^{j}_{m}(\pi,y',s')-\Psi^{j}_{m}(\tilde \pi,y',s')\right|dy'ds'\\
&\leq&2\sum_{m=0}^{q-1}\sum_{i=m+1}^{q}\int_{\ts_{m}}^{\ts_{m+1}}\sum_{j=1}^{q}\int_{\R^{d}}|\pi^{i}-\tilde\pi^{i}|\lambda(\Phi(x_i,s'))e^{-\Lambda(x_i,s')}\\
&&\times Q(\Phi(x_i,s'),x_j)f_{W}(y'-\varphi(x_{j}))dy'ds'\\
&\leq&2\sum_{m=0}^{q-1}\sum_{i=m+1}^{q}|\pi^{i}-\tilde\pi^{i}|\int_{\ts_{m}}^{\ts_{m+1}}\lambda(\Phi(x_i,s'))e^{-\Lambda(x_i,s')}ds'\\
&\leq&2\sum_{i=1}^{q}|\pi^{i}-\tilde\pi^{i}|\int_{0}^{\ts_{i}}\lambda(\Phi(x_i,s'))e^{-\Lambda(x_i,s')}ds'.
\end{eqnarray*}
We obtain the result as $\int_{0}^{\ts_{i}}\lambda(\Phi(x_i,s'))e^{-\Lambda(x_i,s')}ds'=1-e^{-\Lambda(x_{i},\ts_{i})}\leq 1$.
\findemo
\begin{proposition}\label{prop-lip-I}For $v\in BL(\PR(E_{0}))$ and $\left((\pi,u), (\tilde \pi,\tilde u)\right)\in (\PR(E_{0})\times \R^{+})^{2}$, one has
$$|Iv(\pi,u)-Iv(\tilde \pi,\tilde u)|\leq (C_{v}+2[v])|\pi-\tilde \pi|+C_{v}C_{\lambda}|u-\tilde u|.$$
\end{proposition}
\demo 
On the one hand, one clearly has  
$$\left|Iv(\pi,u)-Iv(\pi,\tilde u)\right|\leq \sum_{i=1}^{q}\pi^{i} \big|u\wedge \ts_{i}-\tilde u\wedge \ts_{i}\big|C_{v}C_{\lambda}\leq C_{v}C_{\lambda}|u-\tilde u|.$$
On the other hand, one has
\begin{eqnarray*}\lefteqn{\left|Iv(\pi,u)-Iv(\tilde \pi,u)\right|}\\
&\leq& C_{v}|\pi-\tilde \pi|+\sum_{i=1}^{q}\pi^{i}\int_{0}^{\te_{i}}\int_{\R^{d}}\big|v\big(\Psi(\pi,y',s')\big)-v\big(\Psi(\tilde \pi,y',s')\big)\big|\\
&&\times\sum_{j=1}^{q}Q\big(\Phi(x_{i},s'),x_{j}\big)f_{W}(y'-\varphi(x_{j}))\lambda\circ\Phi(x_{i},s')e^{-\Lambda(x_{i},s')}dy'ds'.
\end{eqnarray*}
Besides, we have assumed that $v$ is Lipschitz continuous so that one has
$$\big|v\big(\Psi(\pi,y',s')\big)-v\big(\Psi(\tilde \pi,y',s')\big)\big|\leq [v]\big|\Psi(\pi,y',s')-\Psi(\tilde \pi,y',s')\big|.$$
Thus, one has
\begin{eqnarray*}\lefteqn{\left|Iv(\pi,y,s,u)-Iv(\tilde \pi,y,s,u)\right|}\\
&\leq& C_{v}|\pi-\tilde \pi|+[v]\sum_{i=1}^{q}\pi^{i}\int_{0}^{\te_{i}}\int_{\R^{d}}\Big|\Psi(\pi,y',s')-\Psi(\tilde \pi,y',s')\Big|\\
&&\sum_{j=1}^{q}Q\big(\Phi(x_{i},s'),x_{j}\big)f_{W}(y'-\varphi(x_{j}))\lambda\circ\Phi(x_{i},s')e^{-\Lambda(x_{i},s')}dy'ds'\\
&\leq& C_{v}|\pi-\tilde \pi|+[v]\sum_{m=0}^{q-1}\sum_{i=m+1}^{q}\pi^{i}\int_{\ts_{m}}^{\ts_{m+1}}\int_{\R^{d}}\Big|\Psi(\pi,y',s')-\Psi(\tilde \pi,y',s')\Big|\\
&&\times\sum_{j=1}^{q}Q\big(\Phi(x_{i},s'),x_{j}\big)f_{W}(y'-\varphi(x_{j}))\lambda\circ\Phi(x_{i},s')e^{-\Lambda(x_{i},s')}dy'ds'\\
&\leq& C_{v}|\pi-\tilde \pi|+[v]\sum_{m=0}^{q-1}\int_{\ts_{m}}^{\ts_{m+1}}\int_{\R^{d}}\Big|\Psi(\pi,y',s')-\Psi(\tilde \pi,y',s')\Big|\overline{\Psi}_{m}(\pi,y',s')dy'ds'.
\end{eqnarray*}
The previous lemma provides the result.
\findemo
\begin{proposition}\label{prop-lip-G}
For $m\in M$, $v\in BL(\PR(E_{0}))$ and $\left((\pi,u), (\tilde \pi,\tilde u)\right)\in (\PR(E_{0})\times \R^{+})^{2}$, one has
$$|G^{m}v(\pi,u)-G^{m}v(\tilde \pi,\tilde u)|\leq (2C_{v}+2[v])|\pi-\tilde \pi|+C_{v}C_{\lambda}|u-\tilde u|.$$
\end{proposition}
\demo 
As in the proof of Proposition \ref{prop-lip-H}, we may assume without loss of generality that $u$, $\tilde u\in [\ts_{m};\ts_{m+1}]$ so that one has
\begin{multline*}
G^{m}v(\pi,u)=Iv(\pi,u)+ \sum_{i=1}^{m}\pi^{i}e^{-\Lambda(x_{i},\te_{i})}\\ 
\times\sum_{j=1}^{q}\int_{\R^{d}}v\big(\Psi(\pi,y',\ts_{i})\big)Q\big(\Phi(x_{i},\te_{i}),x_{j}\big)f_{W}(y'-\varphi(x_{j}))dy',
\end{multline*}
and similarly for $G^{m}v(\tilde \pi,\tilde u)$. The second term does not depend on $u$ thus
\begin{eqnarray*}
\left|G^{m}v(\pi,u)-G^{m}v(\pi,\tilde u)\right|
&=& \left|Iv(\pi,u)-Iv(\pi,\tilde u)\right|\\
&\leq &\left|Iv(\pi,u)-Iv(\tilde \pi,u)\right|+ C_{v}|\pi-\tilde \pi|,
\end{eqnarray*}
as $\Psi(\pi,y',\ts_{i})=\Psi(\tilde \pi,y',\ts_{i})$ by Proposition \ref{prop-rec-filtre}. This yields the result.
\findemo
\begin{proposition}\label{prop-lip-K}
For all $v\in BL(\PR(E_{0}))$ and $(\pi,\tilde \pi) \in \PR(E_{0})^{2}$, one has
$$\left|Kv(\pi)-Kv(\tilde \pi)\right|\leq (2C_{v}+2[v])|\pi-\tilde \pi|.$$
\end{proposition}
\demo 
As $Kv(\pi)=Gv(\pi,\ts_{q})$, this is a consequence of Proposition \ref{prop-lip-G}.
\findemo
\begin{proposition}\label{prop-lip-L}
For $v\in BL(\PR(E_{0}))$ and $(\pi,\tilde \pi) \in \PR(E_{0})^{2}$, one has
$$\left|L(v,g)(\pi)-L(v,g)(\tilde \pi)\right|\leq \left(C_{g}+2C_{v}+2[v]\right) |\pi-\tilde \pi|.$$
\end{proposition}
\demo 
One has
\begin{eqnarray*}
\lefteqn{|L(v,g)(\pi)-L(v,g)(\tilde \pi)|}\\
&\leq&\max_{m\in M}\big\{\sup_{u\in [\ts_{m};\ts_{m+1}[} |J^{m}(v,g)(\pi,u)-J^{m}(v,g)(\tilde \pi,u)|\big\}\vee |Kv(\pi)-Kv(\tilde \pi)|\\
&\leq& \left(C_{g}+2C_{v}+2[v]\right)|\pi-\tilde \pi|,
\end{eqnarray*}
using Propositions \ref{prop-lip-H}, \ref{prop-lip-G} and \ref{prop-lip-K} since $J^{m}(v,g)=H^{m}g+G^{m}v$.
\findemo
\begin{proposition}\label{prop-lip-v}
For all $n\in\{0,\ldots,N\}$, one has $v_{n}\in BL(\PR(E_{0}))$ with
$C_{v_{n}}\leq C_{g}$ and
$[v_{n}]\leq (2^{N-n+2}-3)C_{g}$.
\end{proposition}
\demo
We proved that $v_{n}$ is the value function of the optimal stopping problem with horizon $T_{N-n}$ thus one has
$v_{n}(\pi)=\sup_{\sigma\in \Sigma^{Y}_{N-n}} \E[g(X_{\sigma})\big|\Pi_{0}=\pi]\leq C_{g}.$
Therefore $v_{n}$ is bounded and $C_{v_{n}}\leq C_{g}$.
The second assessment is proved by backward induction. Let $\pi$, $\tilde \pi\in\PR(E_{0})$. One has 
$$|v_{N}(\pi)-v_{N}(\tilde \pi)|\leq \sum_{j=1}^{N}g(x_{j})|\pi^{j}-\tilde \pi^{j}|\leq C_{g}|\pi-\tilde \pi|.$$
Therefore, we have the result for $n=N$ with $[v_{N}]\leq C_{g}$. Moreover, since $v_{n}=L(v_{n+1},g)$ for $0\leq n\leq N-1$, Proposition~\ref{prop-lip-L} yields $[v_{n}]\leq 3C_{g}+2[v_{n+1}]$ which proves the propagation of the induction. 
\findemo


\begin{thebibliography}{10}

\bibitem{arjas92}
E.~Arjas, P.~Haara, and I.~Norros.
\newblock Filtering the histories of a partially observed marked point process.
\newblock {\em Stochastic Process. Appl.}, 40(2):225--250, 1992.

\bibitem{bally03}
V.~Bally and G.~Pag{\`e}s.
\newblock A quantization algorithm for solving multi-dimensional discrete-time
  optimal stopping problems.
\newblock {\em Bernoulli}, 9(6):1003--1049, 2003.

\bibitem{bally05}
V.~Bally, G.~Pag{\`e}s, and J.~Printems.
\newblock A quantization tree method for pricing and hedging multidimensional
  {A}merican options.
\newblock {\em Math. Finance}, 15(1):119--168, 2005.

\bibitem{bauerle11}
N.~B{\"a}uerle and U.~Rieder.
\newblock {\em Markov decision processes with applications to finance}.
\newblock Universitext. Springer, Heidelberg, 2011.

\bibitem{brandejsky11}
A.~Brandejsky, B.~de~Saporta, and F.~Dufour.
\newblock Numerical method for expectations of piecewise-deterministic markov
  processes.
\newblock {\em CAMCoS}, 7(1):63--104, 2012.

\bibitem{bremaud81}
P.~Br{\'e}maud.
\newblock {\em Point processes and queues}.
\newblock Springer Series in Statistics. Springer-Verlag, New York, 1981.

\bibitem{chafai10}
D.~Chafa{\"{\i}}, F.~Malrieu, and K.~Paroux.
\newblock On the long time behavior of the {TCP} window size process.
\newblock {\em Stochastic Process. Appl.}, 120(8):1518--1534, 2010.

\bibitem{davis93}
M.~H.~A. Davis.
\newblock {\em Markov models and optimization}, volume~49 of {\em Monographs on
  Statistics and Applied Probability}.
\newblock Chapman \& Hall, London, 1993.

\bibitem{saporta10}
B.~de~Saporta, F.~Dufour, and K.~Gonzalez.
\newblock Numerical method for optimal stopping of piecewise deterministic
  {M}arkov processes.
\newblock {\em Ann. Appl. Probab.}, 20(5):1607--1637, 2010.

\bibitem{JRR12}
B.~de~{S}aporta, F.~{D}ufour, H.~{Z}hang, and C.~{E}legbede.
\newblock {O}ptimal stopping for the predictive maintenance of a structure
  subject to corrosion.
\newblock {\em Journal of Risk and Reliability}, 226(2):169--181, 2012.

\bibitem{edoli10}
E.~Edoli and W.~J. Runggaldier.
\newblock On optimal investment in a reinsurance context with a point process
  market model.
\newblock {\em Insurance Math. Econom.}, 47(3):315--326, 2010.

\bibitem{gray98}
R.~M. Gray and D.~L. Neuhoff.
\newblock Quantization.
\newblock {\em IEEE Trans. Inform. Theory}, 44(6):2325--2383, 1998.
\newblock Information theory: 1948--1998.

\bibitem{gugerli86}
U.~S. Gugerli.
\newblock Optimal stopping of a piecewise-deterministic {M}arkov process.
\newblock {\em Stochastics}, 19(4):221--236, 1986.

\bibitem{ludkovski12}
M.~Ludkovski and S.~O. Sezer.
\newblock Finite horizon decision timing with partially observable {P}oisson
  processes.
\newblock {\em Stoch. Models}, 28(2):207--247, 2012.

\bibitem{ouvrard04}
J.-Y. Ouvrard.
\newblock {\em Probabilit{\'e}s 2}.
\newblock Enseignement des math{\'e}matiques. Cassini, 2004.

\bibitem{pages04}
G.~Pag{\`e}s, H.~Pham, and J.~Printems.
\newblock Optimal quantization methods and applications to numerical problems
  in finance.
\newblock In {\em Handbook of computational and numerical methods in finance},
  pages 253--297. Birkh\"auser Boston, Boston, MA, 2004.

\bibitem{pakdaman10}
K.~Pakdaman, M.~Thieullen, and G.~Wainrib.
\newblock Fluid limit theorems for stochastic hybrid systems with application
  to neuron models.
\newblock {\em Adv. in Appl. Probab.}, 42(3):761--794, 2010.

\bibitem{pham05}
H.~Pham, W.~Runggaldier, and A.~Sellami.
\newblock Approximation by quantization of the filter process and applications
  to optimal stopping problems under partial observation.
\newblock {\em Monte Carlo Methods Appl.}, 11(1):57--81, 2005.

\end{thebibliography}
\end{document}